\documentclass[11pt]{article} % or book、report

\usepackage{amsthm,amsfonts,amssymb,amsmath,color,mathtools}

\usepackage{hyperref} % Hyperlink ensures that it can be jumped to

\usepackage[]{changes}

\usepackage{todonotes}

\numberwithin{equation}{section}

\newcommand\s{\sigma}

\newcommand\R{\mathbb R}
 
\def\g{\gamma}
\def\de{\delta}

\def\l{\lambda}

\def\epsilon{\varepsilon}
\def\e{\varepsilon}

\newcommand\br{\begin{rem}}
\newcommand\er{\end{rem}}
\newcommand\bp{\begin{pmatrix}}
\newcommand\ep{\end{pmatrix}}
\newcommand\be{\begin{equation}}
\newcommand\ee{\end{equation}}
\newcommand\ba{\begin{equation}\begin{aligned}}
\newcommand\ea{\end{aligned}\end{equation}}

\newcommand\nn{\nonumber}

\usepackage{layout}

\newcommand{\lna}{\langle \nabla \rangle}

\newtheorem{defi}{Definition}[section]
\newtheorem{theorem}[defi]{Theorem}
\newtheorem{proposition}[defi]{Proposition}
\newtheorem{lemma}[defi]{Lemma}
\newtheorem{corollary}[defi]{Corollary}
\newtheorem{remark}[defi]{Remark}

\numberwithin{equation}{section}

\newcommand{\Line}[2]{\underline{\hbox to #1 cm{#2}}}

\newcommand{\pa}{\partial}

\newcommand{\xwh}{\widehat}
\newcommand{\xwt}{\widetilde}
\newcommand{\xwo}{\overline}
\newcommand{\xwl}{\overline}

\newcommand{\Ff}{\mathcal{F}} % Fourier trans
\renewcommand{\Re}{\mathrm{Re}}
\renewcommand{\Im}{\mathrm{Im}}

\newcommand{\la}[1]{\langle #1\rangle}

\newcommand{\D}[1]{\,\mathrm{d}#1}

\begin{document}
\renewcommand{\abstractname}{Abstract}
\renewcommand{\proofname}{Proof}
\renewcommand{\refname}{References}

\title{Uniform-in-time error estimates for the Klein--Gordon--Schr\"{o}dinger system in the nonrelativistic limit regime}

\author{Yong Lu \footnote{School of Mathematics, Nanjing University, 210093 Nanjing, China;  luyong@nju.edu.cn.}\and
    Zhiwei Zheng \footnote{School of Mathematical Sciences, Zhejiang Normal University, 321004 Jinhua, China; zhengzhiwei@zjnu.edu.cn.}
}

\date{}
\maketitle

\begin{abstract}
    In this paper, we study the Klein--Gordon--Schr\"odinger system in the nonrelativistic regime $\e\to 0$ on $\R^d$ with $d\ge4$, where $\e$ is a small parameter proportional to  the inverse of the speed of light. For sufficiently regular initial data, we show that its solution is well approximated by oscillatory profiles governed by a decoupled linear Schr\"odinger-wave semi-limiting model with error $O(\e^2)$ uniformly for all $t \ge 0 $. This rigorously justifies the numerical study of Bao and Zhao \cite{BZ17}.
\end{abstract}

Keywords: Klein--Gordon--Schr\"odinger system, nonrelativistic limit,  Schr\"odinger-wave semi-limiting model, uniform convergence rates.

\section{Introduction}

\subsection{Setting}

In this paper, we consider the dimensionless Klein--Gordon--Schr\"odinger system
\begin{align} \label{4-kgs}
    \begin{dcases}
        \e^2 \pa_t^2 \phi(t,x) - \Delta \phi(t,x) + \frac{1}{\e^2} \phi(t,x)  - \l |\psi(t,x)|^2 = 0, \\
        i \partial_t \psi(t,x) + \Delta \psi(t,x) + \l \phi(t,x) \psi(t,x) = 0.
    \end{dcases}
\end{align}
Here $(t,x)\in\R_+\times\R^d$, $0<\e\ll1$ is proportional to the inverse of the speed of light, and $\l\in\R$ is a coupling constant independent of $\e$. The complex field $\psi$ describes the Schr\"odinger component, while $\phi$ is a real Klein--Gordon field. The coupling models the Yukawa interaction, see Yukawa~\cite{Y35} or Fukuda and Tsutsumi~\cite{FT75,FT78}. For simplicity, we take the mass parameter in \cite{BZ17} to be $\mu=1$. The case of any fixed $\mu>0$ can be treated similarly, with constants depending on $\mu$.

Following the usual nonrelativistic scaling \cite{MNO02,MN02,MN03}, we prescribe the initial data
\begin{align} \label{init-data}
    \phi(0,x) = \phi_0(x), \quad \pa_t \phi(0,x) = \frac{1}{\e^2} \phi_1(x), \quad \psi(0,x) = \psi_0(x),
\end{align}
where $\psi_0$ is complex-valued, $\phi_0$ and $\phi_1$ are real-valued, and all three are independent of $\e$.

The scaling also makes the usual energy estimate unsuitable for a bound uniform in $\e$. Indeed, for smooth solutions the conserved energy is
\begin{align*}
    \int_{\R^d}\big[\frac{\e^2}{2}|\pa_t\phi|^2+\frac{1}{2}|\nabla\phi|^2+\frac{1}{2\e^2}|\phi|^2+|\nabla\psi|^2-\l\phi|\psi|^2\big](t,x)\,\D{x},
\end{align*}
whose quadratic Klein--Gordon part at $t=0$ contains $\frac{1}{2\e^2}(\|\phi_0\|_{L^2}^2+\|\phi_1\|_{L^2}^2)$. The rapid oscillations and this $\e^{-2}$ energy scale motivate estimating the difference between the solution and its oscillatory approximation directly.

For fixed $\e$, the Cauchy problem for \eqref{4-kgs} has been studied extensively. Fukuda and Tsutsumi~\cite{FT75,FT78} and Baillon and Chadam~\cite{BC77} obtained early existence results. Hayashi and von Wahl~\cite{HV87} and Wang~\cite{W97} studied global solutions. Colliander, Holmer, and Tzirakis~\cite{CHT08} proved low-regularity global well-posedness in three dimensions, while Pecher~\cite{P04,P12} obtained further results for rough data. For small, spatially localized data in three dimensions, You~\cite{Y25} proved global existence and scattering by the space-time resonance method.

As $\e\to 0$, the Klein--Gordon component oscillates rapidly in time. Related nonrelativistic and singular-limit problems have been considered for the Maxwell--Klein--Gordon system \cite{BMS04,MN03}, the Maxwell--Landau--Lifshitz system \cite{L13}, the Klein--Gordon--Zakharov system \cite{MN08,MN10,LZ25}, and the Hartree equation \cite{CO06}.
For the coupled system \eqref{4-kgs}, Bao and Zhao~\cite{BZ17} developed a uniformly accurate multiscale numerical method and examined two approximating models: a Schr\"odinger model and a Schr\"odinger-wave model, the latter retaining the $\e^2\pa_t^2$ term in the profile equation.
In \cite{BLZ26}, we established optimal convergence rates for the Schr\"odinger approximation on time intervals of order $\e^{-1}$. In the present paper, we study the Schr\"odinger-wave approximation and prove an $O(\e^2)$ error bound uniform for all $t\ge0$ in dimensions $d\ge4$.

Motivated by the nonrelativistic expansions for Klein--Gordon equations \cite{LZ16,LZ17,MNO02,MN02}, we employ the ansatz
\begin{align}
     & \phi(t,x) =  e^{\frac{it}{\e^2}} \phi_s(t,x) + e^{-\frac{it}{\e^2}} \xwo{\phi_s}(t,x)  + r_1(t,x), \label{4-anstz} \\
     & \psi(t,x) = \psi_s(t,x) + r_2(t,x), \label{4-anstz2}
\end{align}
where the bar denotes complex conjugation, $r_1$ is real-valued, and $r_2$ is complex-valued.

Substituting \eqref{4-anstz} and \eqref{4-anstz2} into the first equation of \eqref{4-kgs}, we obtain
\begin{align} \label{expansion-r1}
    \begin{split}
         & \e^2 \pa_t^2 r_1 - \Delta r_1 + \frac{1}{\e^2} r_1 - \l |r_2|^2 - 2\l \Re(\psi_s \xwl{r_2}) -  \l |\psi_s|^2                                                                                                   \\
         & \quad + e^{\frac{it}{\e^2}} \big( \e^2 \pa_t^2 \phi_s + 2i \pa_t \phi_s - \Delta \phi_s \big) + e^{-\frac{it}{\e^2}} \big( \e^2 \pa_t^2 \xwl{\phi_s} - 2i  \pa_t \xwl{\phi_s} - \Delta \xwl{\phi_s} \big) = 0.
    \end{split}
\end{align}
Substituting the ansatz into the second equation of \eqref{4-kgs}  gives
\begin{align} \label{expansion-r2}
    \begin{split}
         & i \pa_t r_2 + \Delta r_2  + 2\l \Re ( e^{\frac{it}{\e^2}} \phi_s ) r_2 + \l r_1 \psi_s + \l r_1 r_2 + 2\l \Re(e^{\frac{it}{\e^2}} \phi_s ) \psi_s  + i \pa_t \psi_s + \Delta \psi_s = 0.
    \end{split}
\end{align}
To match the initial conditions \eqref{init-data}, the profiles and remainders must satisfy
\begin{align} \label{init-relation}
    \begin{dcases}
        2 \Re \phi_s(0) +r_1(0) = \phi_0,                                                               \\
        -\frac{2 \Im ( \phi_s(0))}{\e^2} + 2 \Re(\pa_t \phi_s(0)) + \pa_t r_1(0) = \frac{\phi_1}{\e^2}, \\
        \psi_s(0) + r_2(0) = \psi_0.
    \end{dcases}
\end{align}

To complete the choice of profiles, we let $\phi_s$ cancel the free oscillatory terms in \eqref{expansion-r1} and $\psi_s$ cancel the free Schr\"odinger term in \eqref{expansion-r2}. The coupling terms remain as forcing for the remainder equations. This gives the following so called semi-limiting linear Schr\"odinger-wave system:
\begin{align} \label{4-semi-eq}
    \begin{dcases}
        \e^2 \pa_t^2 \phi_s + 2i \pa_t \phi_s - \Delta \phi_s = 0, \\
        i \pa_t \psi_s + \Delta \psi_s = 0.
    \end{dcases}
\end{align}

To satisfy the conditions \eqref{init-relation}, we prescribe
\begin{align} \label{init-data-z}
    \phi_s(0,x) = \frac{1}{2}\left( \phi_0(x) - i \phi_1(x) \right), \quad \pa_t \phi_s(0,x) =  0, \quad \psi_s(0,x) = \psi_0(x).
\end{align}
Consequently, all initial data for the remainders vanish:
\begin{align*}
    r_1(0) = 0, \quad \pa_t r_1(0)=0, \quad   r_2(0)=0.
\end{align*}

\subsection{Main results} %\label{sec-result}

The following theorem establishes a uniform in time $O(\e^2)$ approximation of \eqref{4-kgs} by the Schr\"odinger-wave profiles.

\begin{theorem}\label{thm:1}
    Let $d\ge4$, $\s\ge\frac{d}{2}-2$, and $\l\in\R$. Assume that the $\e$-independent initial data $\phi_0,\phi_1$ and $\psi_0$ satisfy
    \begin{align}\label{init:0}
        \phi_0,\ \phi_1\in H^{\s+\frac{5}{2}}(\R^d),
        \quad \psi_0\in H^{\s+2}(\R^d); \quad C_0 := \|(\phi_0, \phi_{1})\|_{H^{\s+\frac52}} + \|\psi_0\|_{H^{\s+2}} >0.
    \end{align}

    \begin{itemize}
        \item[(i)] Then for every $\e \in (0,1]$, the system \eqref{4-semi-eq}--\eqref{init-data-z} admits a unique global solution $(\phi_s,\psi_s)$ satisfying
              \begin{align*}
                  \sup_{t\ge 0}
                  \big(
                  \|\phi_s(t)\|_{H^{\s+\frac52}}
                  +\|\pa_t\phi_s(t)\|_{H^{\s+\frac12}}
                  +\|\psi_s(t)\|_{H^{\s+2}}  \big)
                  \le C C_0,
              \end{align*}
              where $C$ depends only on $d,\s,\l$.

        \item[(ii)] Moreover,  there exists $\e_0\in(0,1]$ depending only on $d,\s,\l$ and the norms in \eqref{init:0} such that
              for every $0<\e\le \e_0$, the system
              \eqref{4-kgs}--\eqref{init-data} admits a unique global solution
              $(\phi,\psi)$ with
              $
                  \phi\in C_b(\R_+;H^{\s+\frac{1}{2}}(\R^d)), \
                  \pa_t\phi\in C_b(\R_+;H^{\s-\frac{1}{2}}(\R^d)), \
                  \psi\in C_b(\R_+;H^\s(\R^d))
              $
              such that
              \begin{align} \label{eq:main-1}
                   & \|\phi(t)-e^{\frac{it}{\e^2}}\phi_s(t)
                  - e^{-\frac{it}{\e^2}}\xwo{\phi_s}(t)\|_{H^\s(\R^d)}
                  +\|\psi(t)-\psi_s(t)\|_{H^\s(\R^d)}\le \xwt{C_0} C_0^2 \e^2, \quad \forall \, t>0.
              \end{align}
              where $\xwt{C_0}$ depends only on $d,\s,\l,C_0$. %Both $\xwt{C_0}$ and $\e_0$ can be chosen uniformly when $C_0$ is bounded above.
    \end{itemize}
\end{theorem}

\begin{remark}
    The condition $\s\ge\frac{d}{2}-2$ ensures that the critical-order remainder estimates needed in the bootstrap argument are available.
    For the Klein--Gordon profile $\phi_s$, the Strichartz estimate used here requires a further half derivative (see Lemma \ref{lem:stri}), which accounts for the assumption $\phi_0,\phi_1\in H^{\s+\frac{5}{2}}(\R^d)$ in \eqref{init:0}.
\end{remark}

%
%\begin{remark}
%    The proof also applies to $\e$-dependent initial data satisfying \eqref{init:0}. In particular, suppose that, for some $\a>0$ and $K>0$ independent of $\e$,
%    \begin{align*}
%        \|\phi_0\|_{H^{\s+\frac52}} + \|\phi_1\|_{H^{\s+\frac52}} + \|\psi_0\|_{H^{\s+2}} \le K \e^\a, \qquad 0<\e\le1.
%    \end{align*}
%    Since the initial-data norms are bounded by $K$, both $\xwt{C_0}$ and $\e_0$ can be chosen uniformly for this family. With the profiles chosen as in \eqref{init-data-z}, estimate \eqref{eq:main-1} therefore gives
%    \begin{align*}
%        \sup_{t\ge0} \bigl(\|\phi(t)-e^{\frac{it}{\e^2}}\phi_s(t)-e^{-\frac{it}{\e^2}}\xwo{\phi_s}(t)\|_{H^\s} + \|\psi(t)-\psi_s(t)\|_{H^\s}\bigr) \le C_K \e^{2\a+2}, \qquad 0<\e\le\e_0.
%    \end{align*}
%\end{remark}
%

\begin{remark}
    \begin{itemize}
        \item Not only the uniform-in-time $O(\e^{2})$ error estimate is given in \eqref{eq:main-1},  but also the clear dependence on norms of initial data is obtained; particularly $C_{0}^{2}$ is kept in the derivation.

        \item   We can also include the coupling strength $|\l|$ in the error estimates.  For $\l\ne 0$, setting $\xwt{\phi}=\l \phi$ and $\xwt{\psi}=\l\psi$ reduces \eqref{4-kgs}--\eqref{init-data} to the system
              \begin{align*}
                  \begin{dcases}
                      \e^2\pa_t^2\xwt{\phi}-\Delta\xwt{\phi}+\e^{-2}\xwt{\phi}-|\xwt{\psi}|^2=0, \\
                      i\pa_t\xwt{\psi}+\Delta\xwt{\psi}+\xwt{\phi}\xwt{\psi}=0,
                  \end{dcases}
              \end{align*}
              with the initial-data $(\l\phi_0,\l\phi_1,\l \psi_0)$.
              Applying Theorem~\ref{thm:1} with coupling constant $1$ and initial-data norm $|\l|C_0$, and scaling back to the original variables, we obtain
              \begin{align*}
                  \sup_{t\ge 0} \bigl(\|\phi(t)-e^{\frac{it}{\e^2}}\phi_s(t)-e^{-\frac{it}{\e^2}}\xwo{\phi_s}(t)\|_{H^\s}+\|\psi(t)-\psi_s(t)\|_{H^\s}\bigr)\le C_\l |\l| \e^2,\qquad 0< \e \le \e_0,
              \end{align*}
              where $C_\l$ and $\e_0$ depend only on $d,\s,|\l|, C_0$. Thus the coupling strength $|\l|$ appears as an additional factor in the $O(\e^2)$ error bound. This is consistent with the finite-time temporal error bound $O(\e^p\tau^2)$ obtained by Bao, Feng, and Su~\cite{BFS22} for the Klein--Gordon equation with $O(\e^p)$ nonlinearity strength.

    \end{itemize}

\end{remark}

Bao and Zhao \cite{BZ17} numerically observed that the Schr\"odinger-wave approximation converges at order $O(\e^2)$, uniformly for all $t\ge 0$, for smooth, well-prepared data
\begin{align*}
    \phi_s(0,x) = \frac{1}{2}\left( \phi_0(x) - i \phi_1(x) \right), \quad  \pa_t \phi_s(0)=-\frac i2\Delta\phi_s(0).
\end{align*}
Here we choose $\pa_t\phi_s(0)=0$ and establish the same uniform error bound rigorously.
A natural question is whether this $O(\e^2)$ convergence rate depends on their particular choice of initial profile velocity. The following corollary shows that it does not: the same error estimate remains valid for any sufficiently regular auxiliary velocity $\Phi_1 \in H^\sigma(\R^d)$.
\begin{corollary}\label{cor:general-Phi1}
    Under the assumptions of Theorem~\ref{thm:1}, let $\Phi_1\in H^\sigma(\R^d)$ be an arbitrary function independent of $\e$, and let
    $(\phi_s^{\Phi_1},\psi_s)$ be the solution of \eqref{4-semi-eq} with
    \begin{align*}
        \phi_s^{\Phi_1}(0)
        =\frac12(\phi_0-i\phi_1),
        \qquad
        \pa_t\phi_s^{\Phi_1}(0)=\Phi_1,
        \qquad
        \psi_s(0)=\psi_0.
    \end{align*}
    Then the corresponding Schr\"odinger-wave approximation satisfies, for every $0< \e \le \e_0$ and $t\ge 0$,
    \begin{align*}
        \|
        \phi(t)- e^{\frac{it}{\e^2}}{\phi_s^{\Phi_1}}(t) - e^{-\frac{it}{\e^2}}\xwo{\phi_s^{\Phi_1}}(t)
        \|_{H^\sigma}
        +  \| \psi(t)-\psi_s(t)
        \|_{H^\sigma}
        \le C_{\Phi_{1}}\, \e^2,
    \end{align*}
    where $C_{\Phi_{1}}$ may depend on $\|\Phi_1\|_{H^\sigma}$ but is independent of
    $t$ and $\e$.
\end{corollary}

For the cubic nonlinear Klein--Gordon equation, the Schr\"odinger approximation was justified in \cite{BLZ24} with an error bound of order $(1+t)\e^2$ up to times of order $\e^{-1}$ in the defocusing case with smooth data.
Lei and Wu~\cite{LW23} obtained a time-uniform $O(\e^2)$ bound for Schr\"odinger-wave profiles in two dimensions and also studied Schr\"odinger profiles in two and three dimensions. Shen, Wang, and Zheng~\cite{SWZ26} constructed higher-order expansions for the nonlinear Klein--Gordon equation. The Strichartz estimates and intervalwise bootstrap arguments in the latter two works are closely related to those used below.

For the Klein--Gordon--Schr\"odinger system, Bao, Dong, and Wang~\cite{BDW10} studied singular limits to Schr\"odinger--Yukawa equations and obtained convergence rates for well-prepared data in dimensions one to three.
Shi and Wang~\cite{SW18} established energy-space convergence for damped and conservative Klein--Gordon--Schr\"odinger systems.
Their scalings and limiting models differ from those considered here. Compared with the numerical study in \cite{BZ17}, Theorem~\ref{thm:1} provides a rigorous analysis on $\R^d$ for $d\ge 4$, with an error bound uniform in $t\ge 0$.

% The rest of the paper is organized as follows. Subsection~\ref{sec:preliminaries} introduces the function spaces and collects the analytic tools used below. Section 2 proves Theorem~\ref{thm:1}(i), establishing global solvability and uniform estimates for the approximate system. Section 3 proves Theorem~\ref{thm:1}(ii) by estimating the remainders, which also yields global existence and uniqueness for the original system.

Throughout the paper, $C$ denotes a positive constant depending only on the dimension, the regularity and integrability exponents, and the fixed parameter $\l$. It is independent of the initial data, $\e$ and time, and may change from line to line unless explicitly fixed. The quantity $C_0$ is the initial-data norm defined in Theorem~\ref{thm:1} and remains fixed throughout the paper. We sometimes use $u(t)$ to denote $u(t,\cdot)$.

The rest of the paper is organized as follows. We give the ideas of proof in the coming Section \ref{sec:main-ideas}, and in Section \ref{sec:preliminaries} we introduce the function spaces and collects the analytic estimates used throughout the proof. Section 2 establishes global solvability and uniform estimates for the approximate system and proves Theorem~\ref{thm:1}(i). Section 3 establishes global existence and uniqueness for the original system and proves the uniform error estimate in Theorem~\ref{thm:1}(ii).

\subsection{Main ideas of the proof} \label{sec:main-ideas}

We prove Theorem~\ref{thm:1} by studying the stability of the Schr\"odinger-wave approximation. Substituting the decompositions \eqref{4-anstz}--\eqref{4-anstz2} into \eqref{4-kgs} gives the remainder equations
\begin{align*}
    \begin{dcases}
        \e^2 \pa_t^2 r_1 - \Delta r_1 + \frac{1}{\e^2} r_1 = \l |r_2|^2 + 2\l \Re(\psi_s \xwl{r_2}) + \l |\psi_s|^2,                                 \\
        i\pa_t r_2 + \Delta r_2 = -2\l \Re(e^{\frac{it}{\e^2}}\phi_s) r_2 - \l r_1 \psi_s - \l r_1 r_2 - 2\l \Re(e^{\frac{it}{\e^2}} \phi_s) \psi_s, \\
        r_1(0) = \pa_t r_1(0) = r_2(0) = 0.
    \end{dcases}
\end{align*}
The first difficulty is to obtain an $O(C_0^2 \e^2)$ bound the pure forcing terms
\begin{align*}
    \l |\psi_s|^2, \qquad -2\l \Re(e^{it/\e^2}\phi_s)\psi_s.
\end{align*}
Both forcing terms are quadratic in the profiles, which accounts for the factor $C_0^2$. To obtain the additional factor $\e^2$, we exploit the rapid temporal oscillations in the Duhamel integrals.

For the Klein--Gordon forcing $\l |\psi_s|^2$, integration by parts in time produces the operator $\e^2(1-\e^2\Delta)^{-1}$. Since $(1-\e^2\Delta)^{-1}$ is uniformly bounded on $H^\s$, we obtain the desired gain of $\e^2$.
For the Schr\"odinger forcing, the corresponding integration by parts is more delicate. Consider the Duhamel term
\begin{align*}
    \int_0^t e^{i(t-s)\Delta}\Re\bigl(e^{\frac{is}{\e^2}}\phi_s(s)\bigr)\psi_s(s) \D{s}.
\end{align*}
For the component with phase $e^{is/\e^2}$, the inverse phase operator $(\Delta-\e^{-2})^{-1} = -\e^2(1-\e^2\Delta)^{-1}$ is uniformly bounded and provides the factor $\e^2$ without a frequency decomposition
However,
\begin{align*}
    \int_0^t e^{i(t-s)\Delta}e^{-\frac{is}{\e^2}}(\xwl{\phi_s}\psi_s)(s) \D{s} & = i\left[e^{i(t-s)\Delta}e^{-\frac{is}{\e^2}}(\Delta+\e^{-2})^{-1}(\xwl{\phi_s}\psi_s)(s)\right]_{s=0}^{s=t}    \\
                                                                               & \qquad-i\int_0^t e^{i(t-s)\Delta}e^{- \frac{i s}{\e^2}}(\Delta+\e^{-2})^{-1}\pa_s(\xwl{\phi_s}\psi_s)(s) \D{s}.
\end{align*}
The identity $(\Delta+\e^{-2})^{-1} = \e^2(1+\e^2\Delta)^{-1}$ suggests the same gain of $\e^2$, but $(1+\e^2\Delta)^{-1}$ is singular at $|\xi| = \e^{-1}$. We therefore separate the low and high spatial frequencies.

At low frequencies $|\xi|\le\frac{1}{2\e}$, the inverse phase factors are bounded by $C\e^2$. At high frequencies, we estimate the Duhamel term directly and use two additional spatial derivatives to obtain the factor $\e^2$. For a smooth high-frequency projection $P_{>\frac{1}{2\e}}$, we have
\begin{align*}
    \|P_{>\frac{1}{2\e}}(\xwl{\phi_s}\psi_s)\|_{H^\s} & \le C \e^2 \||\nabla|^2(\xwl{\phi_s}\psi_s)\|_{H^\s} \le C \e^2 \|\xwl{\phi_s}\psi_s\|_{H^{\s+2}}.
\end{align*}
Combining these frequency-localized bounds with the spacetime product estimates and Strichartz estimates, we obtain a pure forcing contribution of order $C_0^2\e^2$.

The second difficulty is to control the coupling terms uniformly for all time.
Although the profiles $\phi_s$ and $\psi_s$ have finite global critical spacetime norms, these norms need not be small. Thus the mixed profile-remainder terms
\begin{align*}
    2\Re(e^{it/\e^2}\phi_s) r_2, \qquad r_1 \psi_s, \qquad 2\Re(\psi_s \xwl{r_2}),
\end{align*}
cannot  be absorbed by applying the estimates directly on $\R_+$. We thus divide $\R_+$ into finitely many intervals
\begin{align*}
    0 = t_0 < t_1 < \cdots < t_n = \infty, \qquad I_k = [t_k,t_{k+1}), \quad k = 0,1,\ldots,n-1,
\end{align*}
on which the critical spacetime norms of the profiles are sufficiently small. At the critical order $s_c = \frac d2-2$, this smallness allows us to absorb the mixed profile-remainder terms. The remaining quadratic terms $|r_2|^2$ and $r_1 r_2$ are controlled by a bootstrap argument using the bound $C C_0^2 \e^2$ for the pure forcing contributions. We then estimate the orders $0$ and $\s$ and propagate the endpoint bounds from one interval to the next.

Since the number of intervals has an upper bound independent of $\e$, we can choose $\e_0$ sufficiently small so that the bootstrap closes on every interval. The continuation criterion then yields a global solution satisfying
\begin{align*}
    \sup_{t\ge 0} \bigl(\|r_1(t)\|_{H^\s}+\|r_2(t)\|_{H^\s}\bigr) \le \xwt{C_0} C_0^2 \e^2.
\end{align*}

\subsection{Preliminaries}\label{sec:preliminaries}

In this subsection, we recall some notation, Strichartz estimates and the fractional Leibniz rule, which will be frequently used in the arguments below.

For $u$ on $\R^d$, we write its Fourier transform as
\begin{align*}
    \xwh{u} (\xi) =\Ff u (\xi) = \int_{\R^d} e^{-i x\cdot\xi} u(x) \D{x}.
\end{align*}
Let $\chi$ be a smooth radial cutoff satisfying
\begin{align*}
    \chi(\xi)=1 \quad\text{if } |\xi|\leq \frac{1}{2},\quad  \operatorname{supp}\chi \subset \{\xi:|\xi|\leq 1 \},\quad
    0\leq \chi(\xi)\leq 1.
\end{align*}
For $K>0$, define $P_{\leq K}$ by
\begin{align*}
    \widehat{P_{\leq K}u}(\xi) := \chi(\frac{\xi}{K} ) \hat u(\xi),
\end{align*}
and set $P_{>K}u := u-P_{\leq K}u$.

We recall the definition of a Fourier multiplier $a(D)$ related to a symbol $a$:
$$
    \widehat{a(D) u}(\xi)  : = a(\xi) \hat u (\xi).
$$
In particular, $|\nabla|, \ \Delta$ and  $\lna: = (1-\Delta)^{\frac 12} =(1 + |\nabla|^{2})^{\frac 12}$ are Fourier multipliers with symbols $|\xi|,  \  -|\xi|^{2},$ and $(1+|\xi|^{2})^{\frac 12}$, respectively.

We first recall the Mikhlin multiplier theorem, see, for instance, \cite{S70,G14}.
\begin{lemma}
    Let $1<p<\infty$, and let $a \in C^N(\R^d\setminus\{0\})$ for some integer $N > \frac{d}{2}$. Suppose that
    \begin{align*}
        \sup_{\xi\neq 0} |\xi|^{|\alpha|}
        \left|\partial_\xi^\alpha a(\xi)\right|
        \le C_\alpha,
    \end{align*}
    for every multi-index \(\alpha\) with \(|\alpha|\le N\). Then the Fourier multiplier operator $a(D)$
    extends to a bounded operator on $L^p(\mathbb \R^d)$. In particular,
    \begin{align*}
        \|a(D) f\|_{L^p} \le C\|f\|_{L^p},
    \end{align*}
    where $C$ depends only on $d,p,N$ and $\max_{|\alpha|\le N}C_\alpha$.
\end{lemma}
We also recall the fractional Leibniz rule, also known as the Kato--Ponce inequality \cite{GO14}.
\begin{lemma} \label{lem:fractional-Leibniz}
    Let $s>0$, $1 < p < \infty$, and let $1<p_1,p_2,q_1,q_2\le \infty$ satisfy
    \begin{align*}
        \frac1p=\frac1{p_1}+\frac1{q_1} =\frac1{p_2}+\frac1{q_2}.
    \end{align*}
    Then
    \begin{align*}
        \| |\nabla|^s (fg)\|_{L^p}
        \le C \big(
        \||\nabla|^s f\|_{L^{p_1}}\| g \|_{L^{q_1}}
        +
        \|f\|_{L^{p_2}}\||\nabla|^s g\|_{L^{q_2}} \big),
    \end{align*}
    where the implicit constant depends only on \(d,s,p,p_1,p_2,q_1,q_2\).
\end{lemma}

Next, we introduce the mixed time-space norms and the Strichartz estimate.
For a time interval $I$, we write
\begin{align*}
    \| f \|_{L^p_{t} L^q_{x}(I)} := \Big(\int_I \Big( \int_{\R^d}  |f(t,x)|^q \D{x} \Big)^{\frac{p}{q}} \D{t} \Big)^{\frac{1}{p}}.
\end{align*}

A pair of exponents $(p,q)$ is called Schr\"odinger-admissible if $2\leq p,q\leq \infty$, $(p,q,d)\ne (2,\infty,2)$, and
\begin{align*}
    \frac{2}{p} + \frac{d}{q} = \frac{d}{2}.
\end{align*}

Let $e^{it\Delta}$ be the linear Schr\"odinger group and let $e^{it\lna}$ be the linear Klein--Gordon group. We recall the following Strichartz estimates (see \cite{KT98} for the Schr\"odinger group and \cite{CGM23} for the Klein-Gordon group):
\begin{lemma} \label{lem:stri}
    Let $I=[a,b) \subset \R_+$, allowing $b=\infty$. Suppose that $(p,q)$ is a Schr\"odinger-admissible pair. Then,
    \begin{align}
        (i): \| e^{it\Delta} f \|_{L_t^p L_x^q(I)} & \le C \| f \|_{L_x^2(\R^d)},             \label{est:s-1}                                     \\
        (ii): \| e^{it\lna} f \|_{L_t^p L_x^q(I)}  & \le C \| \lna^{\frac{d+2}{2}(\frac{1}{2}-\frac{1}{q})} f \|_{L_x^2(\R^d)}.  \label{est:kg-1}
    \end{align}
    If, in addition, $(\xwt{p}, \xwt{q})$ is a Schr\"odinger-admissible pair, then the following estimates for the Duhamel terms hold:
    \begin{align}
        (iii): \| \int_a^t e^{i(t-s)\Delta} F(s) \D{s} \|_{L_t^p L_x^q(I)} & \le C \| F \|_{L_t^{\xwt{p}'} L_x^{\xwt{q}'}(I)},  \label{est:s-2}                                                      \\
        (iv):  \| \int_a^t e^{i(t-s) \lna} F(s) \D{s} \|_{L_t^p L_x^q(I)}  & \le C \| \lna^{\frac{d+2}{2}(1-\frac{1}{q}-\frac{1}{\xwt{q}})} F\|_{L_t^{\xwt{p}'} L_x^{\xwt{q}'}(I)}. \label{est:kg-2}
    \end{align}
    Here $\xwt{p}'$ and $\xwt{q}'$ denote the H\"older conjugates of $\xwt{p}$ and $\xwt{q}$, respectively.
\end{lemma}

We now introduce the space-time norms used throughout the paper.
Let $d\geq 4$ and
$$
    p_{1} := \frac{d+2}{2}, \quad q_{1}: = \frac{2d(d+2)}{d^2+2d-8}; \quad p_{2} = q_{2}: = \frac{2(d+2)}{d}.
$$
Clearly $  (p_{1}, q_{1})$ and $(p_{2}, q_{2})$ are both Schr\"odinger-admissible pairs.
For $I\subset\R_+$, define
\begin{align}  \label{def:norm-X}
    \| f \|_{X(I)} := \| f \|_{L_t^{p_{1}} L_x^{q_{1}}(I)} + \|  f \|_{L_t^{p_{2}} L_x^{q_{2}}(I)},
\end{align}
and
\begin{align*}  %\label{def:Y}
    \| f \|_{Y^0(I)} := \| f \|_{L^\infty L^2(I) }+ \| f \|_{X(I)}, \quad \| f \|_{Y^{\frac{1}{2}}(I)} := \| \lna^{\frac{1}{2}} f \|_{L^\infty L^2(I) }+ \| f \|_{X(I)}.
\end{align*}
It follows directly from the definitions that
\begin{align*}
    Y^{\frac12}(I)\subset Y^0(I)\subset X(I).
\end{align*}
We also use the fact that the finite time exponents defining $X$ imply absolute continuity with respect to the time interval.
More precisely, if $f\in X(\R_+)$, then for every $\de > 0$, there exists $\eta>0$ such that
\begin{align*} %\label{lem:X-ac}
    \|f\|_{X(J)}< \de,
\end{align*}
whenever $J \subset\R_+$ is an interval with $|J|<\eta$.

The preceding Strichartz estimates yield the following bounds in the spaces defined above.
\begin{lemma}
    Let $d\ge4$ and let $I=[a,b)\subset\R_+$, allowing $b=\infty$.
    The following estimates hold.
    \begin{itemize}
        \item   For the  Schr\"odinger group,
              \begin{align} \label{est:s-3}
                  \| e^{it\Delta} f \|_{Y^0(I)} \le C \| f \|_{L_x^2(\R^d)}, \quad \| \int_a^t e^{i(t-s)\Delta} F(s) \D{s} \|_{Y^0(I)} \le C \|F \|_{L_{tx}^{\frac{2(d+2)}{d+4}}(I)}.
              \end{align}

        \item    For the Klein--Gordon group,
              \begin{align}
                  \| e^{it\lna} f\|_{Y^{\frac{1}{2}}(I)}                                    & \le C \| \lna^{\frac{1}{2}} f \|_{L_x^2(\R^d) },  \label{est:kg-3}    \\
                  \| \int_a^t e^{i(t-s)\lna } \lna^{-1} F(s) \D{s} \|_{ Y^{\frac{1}{2}}(I)} & \le C \|  F \|_{L_{tx}^{\frac{2(d+2)}{d+4}}(I) }.    \label{est:kg-4}
              \end{align}

    \end{itemize}

\end{lemma}

\begin{proof}
    The assertion for the Schr\"odinger group follows directly from \eqref{est:s-1} and \eqref{est:s-2}.
    For the Klein--Gordon group, the condition $d\ge 4$ gives
    \begin{align*}
        \|\lna^{\frac{2}{d}}f\|_{L_x^2(\R^d) }  \le  \|\lna^{\frac{1}{2}}f\|_{L_x^2(\R^d) }.
    \end{align*}
    Thus, applying \eqref{est:kg-1} to $X$ yields
    \begin{align*}
        \| e^{it\lna} f\|_{Y^{\frac{1}{2}}(I)} \le C \| \lna^{\frac{1}{2}} f \|_{L_x^2(\R^d) } + C  \| \lna^{\frac{2}{d}} f \|_{L_x^2(\R^d) } \le C \| \lna^{\frac{1}{2}} f \|_{L_x^2(\R^d) }.
    \end{align*}
    Since $d\ge 4$, we have $\frac2d-\frac12\le0$. The Mikhlin multiplier theorem therefore implies that
    $\lna^{\frac2d-\frac12}$ is bounded on $L_x^{\frac{2(d+2)}{d+4}}(\R^d)$.
    Applying \eqref{est:kg-2} to the components of $Y^{\frac12}(I)$ yields
    \begin{align*}
        \| \int_a^t e^{i(t-s)\lna } \lna^{-1} F(s) \D{s} \|_{ Y^{\frac{1}{2}}(I)} \le C \|F \|_{L_{tx}^{\frac{2(d+2)}{d+4}}(I)} + C   \|\lna^{\frac{2}{d}-\frac{1}{2}} F \|_{L_{tx}^{\frac{2(d+2)}{d+4}}(I)} \le  C \|F\|_{L_{tx}^{\frac{2(d+2)}{d+4}}(I)}.
    \end{align*}
    % The constant $C$ depends only on $d$ and is independent of the endpoints of $I$.
\end{proof}

We conclude this subsection with two product estimates used in the remainder analysis. The first is a spatial estimate obtained from the fractional Leibniz rule and Sobolev embedding.
\begin{lemma}\label{ineq:eta}
    Let $2 \le q < \infty$, $0<\eta<\frac{d}{q}$ and $\g \ge 0$. Then the following estimate holds for all sufficiently regular functions $(f,g)$:
    \begin{align*}
        \| |\nabla|^\gamma (fg) \|_{L_x^q(\R^d)}
        \le C \| |\nabla|^{\gamma+\eta} f \|_{L_x^q(\R^d)}
        \| |\nabla|^{\frac{d}{2}-\eta} g \|_{L_x^2(\R^d)}
        + C  \| |\nabla|^{\frac{d}{2}-\eta} f \|_{L_x^2(\R^d)}
        \| |\nabla|^{\gamma+\eta} g \|_{L_x^q(\R^d)}.
    \end{align*}
    The constant $C$ depends only on $d,q,\eta,\gamma$.
\end{lemma}

\begin{proof}
    For $\g=0$, the following estimate follows from H\"older's inequality.
    For $\g>0$, it follows from Lemma~\ref{lem:fractional-Leibniz}:
    \begin{align*}
        \||\nabla|^\gamma(fg)\|_{L_x^q(\R^d)}
        \le
        C \||\nabla|^\gamma f\|_{L_x^{\frac{d q}{d-\eta q} }(\R^d)} \|g\|_{L_x^{\frac d\eta }(\R^d)}
        + C
        \|f\|_{L_x^{\frac d\eta }(\R^d)} \||\nabla|^\gamma g\|_{L_x^{ \frac{d q}{d-\eta q} }(\R^d)},
    \end{align*}
    By Sobolev embedding,
    \begin{align*}
        \||\nabla|^\gamma f\|_{L_x^{\frac{d q}{d-\eta q} }(\R^d)}
        \le C \||\nabla|^{\gamma+\eta}f\|_{L_x^q(\R^d)},
    \end{align*}
    while
    \begin{align*}
        \|g\|_{L_x^{\frac{d}{\eta}}(\R^d)}
        \le
        C \||\nabla|^{\frac{d}{2}-\eta}g\|_{L_x^2(\R^d)}.
    \end{align*}
    Applying the same two embeddings with $f$ and $g$ interchanged yields the desired estimate for the second term.

\end{proof}

\begin{remark}
    The restriction $\eta>0$ avoids the endpoint Sobolev embedding in the preceding proof. Setting $\eta=0$ would require $\dot H^{\frac d2}(\R^d) \hookrightarrow L^\infty(\R^d)$, which does not hold in general.
\end{remark}

We use the following bilinear estimate to control the quadratic terms in the Duhamel formulas.
\begin{lemma} \label{lem:2}
    Let $d\ge4$, $s_c=\frac d2-2$, $I\subset\R_+$, and $\g\ge 0$. Then the following estimate holds for all sufficiently regular functions $(f,g)$:
    \begin{align*} %\label{ineq:2d+4}
        \| |\nabla|^\g (fg) \|_{L_{tx}^{\frac{2(d+2)}{d+4}}(I)} \le C \| |\nabla|^\g f \|_{X(I)} \| |\nabla|^{s_c} g\|_{X (I)} + C  \| |\nabla|^{s_c} f \|_{X(I)} \| |\nabla|^\g g \|_{X(I)}.
    \end{align*}
\end{lemma}
\begin{proof}
    For $\g=0$, H\"older's inequality gives the following bound directly. For $\g>0$, the same bound follows from the fractional Leibniz rule and H\"older's inequality:
    \begin{align*}
        \| |\nabla|^\g (fg) \|_{L_{tx}^{\frac{2(d+2)}{d+4}}(I)}
         & \le C \| |\nabla|^\g f \|_{L_{tx}^{\frac{2(d+2)}{d}}(I)  }
        \|  g \|_{L_{tx}^{\frac{d+2}{2}}(I)  }
        + C   \|  f \|_{L_{tx}^{\frac{d+2}{2}} (I)} \| |\nabla|^\g  g \|_{L_{tx}^{\frac{2(d+2)}{d} }(I) }   .
    \end{align*}
    Using the Sobolev embedding
    \begin{align*}
        \dot{W}^{s_c, \frac{2d(d+2)}{d^2+2d-8}}(\mathbb{R}^d) \hookrightarrow L^{\frac{d+2}{2}}(\mathbb{R}^d),
    \end{align*}
    we have
    \begin{align*}
        \|g\|_{L_{tx}^{\frac{d+2}{2}}(I)}
        \le
        C \| |\nabla|^{s_c}g\|_{L_t^{\frac{d+2}{2}}
                L_x^{\frac{2d(d+2)}{d^2+2d-8}}(I) },
    \end{align*}
    and the analogous estimate for $f$. Together with the definition of $X(I)$, this proves the lemma.
\end{proof}

\begin{remark}
    The choice of $(p_{1}, q_{1})$ and $(p_{2}, q_{2})$ exponents may not be optimal in the sense of initial regularity assumption.  For $d\ge4$, replacing both pairs by
    \begin{align*}
        (p,q)=\big(\frac{4(d+2)}{d+4},\frac{2d(d+2)}{d^2+d-4}\big)
    \end{align*}
    gives the Klein--Gordon derivative loss (see \eqref{est:kg-1})
    \begin{align*}
        \frac{d+2}{2}\big(\frac12-\frac1q\big)=\frac14+\frac1d\le\frac12.
    \end{align*}
    With corresponding adjustments to the profile estimates, the regularity assumptions on the Klein--Gordon initial data could be replaced by
    \begin{align*}
        \phi_0,\phi_1\in H^{\s+2+(\frac14+\frac1d)}(\R^d),
    \end{align*}
    which is slightly weaker than our assumption $H^{\s+\frac52}(\R^d)$ for $d>4$. However, we do not pursue this refinement here.
\end{remark}

\section{Estimates for the approximate solution}

In this section, we prove Theorem~\ref{thm:1}(i). We establish global existence, uniqueness, and uniform Sobolev bounds for the profiles solving \eqref{4-semi-eq}--\eqref{init-data-z}. We also derive the space-time estimates needed for the remainder analysis in Section 3.

Define the scaling operator
\begin{align*}
    S_\e f(t,x):=\e^2 f(\e^2t,\e x),
\end{align*}
and, for time-independent functions, $S_\e f(x):=\e^2f(\e x)$. Set
\begin{align} \label{h-f}
    \phi^\e := S_\e \phi_s, \quad \psi^\e := S_\e \psi_s.
\end{align}
Note that for any Schr\"odinger-admissible pair $(p,q)$ and $\g\ge 0$,
\begin{align} \label{scale}
    \| |\nabla|^\g (\phi^\e, \psi^\e) \|_{L_t^p L_x^q (\R_+)} = \e^{\g-s_c} \| |\nabla|^\g (\phi_s, \psi_s) \|_{L_t^p L_x^q(\R_+)}, \quad s_c := \frac{d}{2}-2.
\end{align}
For each fixed $\e>0$, $S_\e$ and its inverse also preserve continuity in $H^\g$ on the corresponding time intervals. %We use this qualitative property when returning to the original variables below.

By \eqref{h-f}, the first equation in \eqref{4-semi-eq} becomes
\begin{align} \label{eq:phi}
    \begin{dcases}
        \pa_t^2 \phi^\e + 2i\pa_t \phi^\e - \Delta \phi^\e =0, \\
        \phi^\e(0,x) = \frac{1}{2} S_\e(\phi_0-i\phi_1),       \\
        \pa_t\phi^\e(0,x) = 0.
    \end{dcases}
\end{align}
The following proposition gives bounds uniform in rescaled time.

\begin{proposition} %\label{prop-Sephi}
    Let $d\ge 4$ and $0<\e\le 1$.  Assume that the initial data satisfy \eqref{init:0}. Then the solution $\phi^\e$ satisfies
    \begin{align}
        (i):  & \quad
        \| |\nabla|^\g \phi^\e \|_{Y^{\frac{1}{2}}(\R_+)}
        \le C C_0 \e^{\g-s_c}, \quad \text{for all $0\le \g \le \s+2$},
        \label{est:phi-0}
        \\
        (ii): & \quad
        \| |\nabla|^\g \pa_t \phi^\e \|_{Y^{\frac{1}{2}}(\R_+)}
        \le  C C_0  \e^{\g+2-s_c},  \quad \text{for all $0\le \g \le \s$}.
        \label{est:phi-1}
    \end{align}
\end{proposition}

\begin{proof}
    Taking the Fourier transform of \eqref{eq:phi} gives
    \begin{align*}
        \pa_t^2 \xwh{\phi^\e} + 2i \pa_t \xwh{\phi^\e}  + |\xi|^2 \xwh{\phi^\e}  =0.
    \end{align*}
    Setting $\xwt{\phi}(t,\xi):=e^{it}\xwh{\phi^\e}(t,\xi)$, we obtain
    \begin{align} \label{eq:phi2}
        \pa_t^2 \xwt{\phi} + (|\xi|^2+1) \xwt{\phi}  =0,
    \end{align}
    with the initial conditions
    \begin{align*} %\label{init:phi2}
        \xwt{\phi}(0,\xi) = \xwh{\phi^\e}(0,\xi), \quad \pa_t \xwt{\phi}(0,\xi) = i \xwh{\phi^\e}(0,\xi) + \pa_t \xwh{\phi^\e}(0,\xi).
    \end{align*}
    Solving \eqref{eq:phi2} with these data gives
    \begin{align*}
        \xwt{\phi}(t,\xi) = \cos(t \la{\xi}) \xwh{\phi^\e}(0,\xi) + \frac{\sin (t \la{\xi})}{\la{\xi}} (i \xwh{\phi^\e}(0,\xi) + \pa_t \xwh{\phi^\e}(0, \xi)).
    \end{align*}
    Inverting the Fourier transform yields
    \begin{align} \label{eq:phie}
        \phi^\e (t,x) = & \frac{1}{2} e^{-it} \big( \cos(t\lna) + i\sin(t\lna) \lna^{-1} \big)  S_\e  (\phi_0-i\phi_1).
    \end{align}
    Applying \eqref{est:kg-3} to the Klein--Gordon propagators in \eqref{eq:phie}, we obtain
    \begin{align*}
        \| |\nabla|^\g \phi^\e  \|_{Y^\frac{1}{2}(\R_+)} & \le C \| \lna^{\frac{1}{2}} |\nabla|^\g S_\e (\phi_0-i\phi_1) \|_{L_x^2(\R^d) }                                                           \\
                                                         & \le  C \|  |\nabla|^\g S_\e (\phi_0-i\phi_1) \|_{L_x^2(\R^d) } +  \|  |\nabla|^{\g+\frac{1}{2}} S_\e (\phi_0-i\phi_1) \|_{L_x^2(\R^d) } .
    \end{align*}
    By the scaling law \eqref{scale},
    \begin{align*}
        \| |\nabla|^\g \phi^\e  \|_{Y^{\frac{1}{2}}(\R_+)} & \le C \e^{\g-s_c} \|   (\phi_0, \phi_1) \|_{H^\g (\R^d) } + C \e^{\g + \frac{1}{2}-s_c} \|   (\phi_0, \phi_1) \|_{H^{\g+\frac{1}{2}}(\R^d) }.
    \end{align*}
    Since $0<\e\le 1$ and $\g+\frac12 \le \s+\frac52$ for $0\le\g\le\s+2$, assumption \eqref{init:0} implies \eqref{est:phi-0}.

    Differentiating \eqref{eq:phie} gives
    \begin{align*}
        \begin{split}
            \pa_t \phi^\e (t,x) & = - \frac{1}{2} e^{-it} \sin(t\lna) \lna^{-1} |\nabla|^2 S_\e (\phi_0-i\phi_1) .
        \end{split}
    \end{align*}
    Since $\lna^{-\frac12}$ is bounded on $L_x^2$ , we obtain, for $0\le \g \le \s$,
    \begin{align*}
        \||\nabla|^\g\pa_t\phi^\e\|_{Y^{\frac12}(\R_+)}
         & \le C\|\lna^{-\frac12}|\nabla|^{\g+2}S_\e(\phi_0-i\phi_1)\|_{L_x^2} \\
         & \le C\||\nabla|^{\g+2}S_\e(\phi_0-i\phi_1)\|_{L_x^2}                \\
         & \le C\e^{\g+2-s_c}\|(\phi_0,\phi_1)\|_{H^{\s+2}}
        \le C C_0\e^{\g+2-s_c}.
    \end{align*}
\end{proof}

The Schr\"odinger profile satisfies
\begin{align} \label{eq:psie}
    \begin{dcases}
        i \pa_t \psi^\e + \Delta \psi^\e = 0, \\
        \psi^\e (0,x) = S_\e \psi_0.
    \end{dcases}
\end{align}
We therefore obtain the corresponding bounds directly from the Schr\"odinger estimates.
\begin{proposition} %\label{prop-Sg}
    Let $d\ge 4$ and $0<\e\le 1$. Assume that the initial data $\psi_0$ satisfy \eqref{init:0}. Let $\psi^\e$ denote the unique global-in-time solution to \eqref{eq:psie}. Then the following estimates hold:
    \begin{align}
        (i):  & \quad  \| |\nabla|^\g   \psi^\e \|_{ Y^0(\R_+) } \le  C C_0 \e^{\g-s_c},  \quad \text{for all $0\leq \g \leq \s +2$}, \label{est:psi-0}     \\
        (ii): & \quad  \| |\nabla|^\g \pa_t \psi^\e \|_{ Y^0(\R_+) } \le C C_0 \e^{\g + 2 -s_c}, \quad \text{for all $0\leq \g \leq \s$}. \label{est:psi-1}
    \end{align}
\end{proposition}
\begin{proof}
    By \eqref{eq:psie},
    \begin{align} \label{formula:psie}
        \psi^\e(t)=e^{it\Delta}S_\e\psi_0.
    \end{align}
    Applying \eqref{est:s-1} to the pairs defining $Y^0$ gives
    \begin{align*}
        \| |\nabla|^\g \psi^\e \|_{Y^0(\R_+)}
        \le
        C \| |\nabla|^\g S_\e\psi_0\|_{L_x^2(\R^d)}.
    \end{align*}
    The scaling identity \eqref{scale} then yields
    \begin{align*}
        \| |\nabla|^\g \psi^\e \|_{Y^0(\R_+)} \le C \e^{\g-s_c} \| |\nabla|^\g \psi_0 \|_{L_x^2(\R^d)} \le  C C_0 \e^{\g-s_c}.
    \end{align*}

    Since $\pa_t\psi^\e=i\Delta\psi^\e$,
    \begin{align*}
        \| |\nabla|^\g\pa_t\psi^\e\|_{Y^0(\R_+)}
        \le
        C \| |\nabla|^{\g+2}S_\e\psi_0\|_{L_x^2(\R^d)}
        \le
        C C_0 \e^{\g+2-s_c},
    \end{align*}
    which proves \eqref{est:psi-1}.
\end{proof}

\begin{proof}[Proof of Theorem~\ref{thm:1}(i)]
    The explicit formulas \eqref{eq:phie} and \eqref{formula:psie}, together with the inverse scaling, give unique global profile solutions with the stated time continuity. To obtain uniform Sobolev bounds in the original variables, we first apply Plancherel's theorem to \eqref{eq:phie} and its time derivative:
    \begin{align*}
        \sup_{t\ge 0} \||\nabla|^\g\phi^\e(t)\|_{L_x^2}
         & \le C\e^{\g-s_c}\|(\phi_0,\phi_1)\|_{H^\g},
        \qquad 0\le\g\le\s+\frac52,                           \\
        \sup_{t\ge 0} \||\nabla|^\g\pa_t\phi^\e(t)\|_{L_x^2}
         & \le C\e^{\g+2-s_c} \|(\phi_0,\phi_1)\|_{H^{\g+2}},
        \qquad 0\le\g\le\s+\frac12.
    \end{align*}
    Applying \eqref{scale} at order zero and at the highest derivative order for each profile, and using $\pa_t\phi^\e(t,x)=\e^2 S_\e(\pa_t\phi_s)(t, x)$ together with \eqref{est:psi-0}, we obtain
    \begin{align*}
        \sup_{t\ge0}\big(
        \|\phi_s(t)\|_{H^{\s+\frac52}}
        +\|\pa_t\phi_s(t)\|_{H^{\s+\frac12}}
        +\|\psi_s(t)\|_{H^{\s+2}}
        \big)
         & \le C\big(
        \|\phi_0\|_{H^{\s+\frac52}}
        +\|\phi_1\|_{H^{\s+\frac52}}
        +\|\psi_0\|_{H^{\s+2}}
        \big)         \\
         & \le C C_0,
    \end{align*}
    uniformly for $0<\e\le1$. This proves part (i) of Theorem~\ref{thm:1}.
\end{proof}

\section{Error estimates}

This section is devoted to the proof of Theorem~\ref{thm:1}(ii).
We first reduce the original problem to the remainder system \eqref{eq-for-R1R2} and establish local existence together with a continuation criterion.
We then construct a finite partition of $\R_+$ on which the profile norms are small, derive a priori estimates for the two remainder components, and close the bootstrap argument.

We now derive the equations for the remainders in \eqref{4-anstz}--\eqref{4-anstz2} and estimate their two components.
Using the profile equations \eqref{4-semi-eq} in \eqref{expansion-r1}--\eqref{expansion-r2}, we obtain
\begin{align} \label{eq:r}
    \begin{dcases}
        \e^2 \pa_t^2 r_1 - \Delta r_1 + \frac{1}{\e^2} r_1 -  \l |r_2|^2 - 2 \l \Re(\psi_s \xwl{r_2}) - \l |\psi_s|^2  = 0, \\
        i \pa_t r_2 + \Delta r_2  + 2\l \Re ( e^{\frac{it}{\e^2}} \phi_s ) r_2 + \l r_1 \psi_s + \l r_1 r_2 + 2\l \Re(e^{\frac{it}{\e^2}} \phi_s ) \psi_s  = 0.
    \end{dcases}
\end{align}
From the initial relation \eqref{init-relation}, we have
\begin{align*} %\label{init:r1r2}
    r_1(0) = 0, \quad \pa_t r_1(0) =  0, \quad  r_2(0) = 0.
\end{align*}

Set
\begin{align*} %\label{tran-Rr-2}
    R_1(t,x) := S_\e r_1(t,x) , \quad R_2(t,x) := S_\e r_2(t,x).
\end{align*}
Using $\phi^\e=S_\e\phi_s$ and $\psi^\e=S_\e\psi_s$, the system \eqref{eq:r} becomes
\begin{align} \label{eq-for-R1R2}
    \begin{dcases}
        \pa_t^2 R_1 - \Delta R_1 + R_1 -  \l |R_2|^2 - 2 \l \Re(\psi^\e \xwl{R_2}) - \l | \psi^\e|^2 = 0, \\
        i \pa_t R_2 + \Delta R_2   + 2\l \Re ( e^{it} \phi^\e ) R_2 + \l R_1  \psi^\e + \l R_1 R_2 + 2\l \Re(e^{it}  \phi^\e )  \psi^\e =0,
    \end{dcases}
\end{align}
with initial data
\begin{align} \label{init-R1R2}
    R_1(0) = 0, \quad \pa_t R_1(0) =  0, \quad R_2(0) = 0.
\end{align}

We introduce the transformation
\begin{align*} %\label{def-R1-x}
    \xwt{R_1} = \lna^{-1} (\pa_t+ i \lna ) R_1.
\end{align*}
Since $R_1$ is real-valued, a direct computation gives
\begin{align} \label{relRr2}
    R_1 = \Im(\xwt{R_1}), \quad \pa_t R_1 = \Re ( \lna \xwt{R_1}).
\end{align}
In particular, for every interval $I\subset\R_+$, $\| R_1 \|_{Y^0(I)} \leq \| \xwt{R_1} \|_{Y^0(I)}$.
Then, the equation for $\xwt{R_1}$ is
\begin{align*}
    \begin{dcases}
        \pa_t\xwt{R_1}-i\lna\xwt{R_1}
        -\lna^{-1}\bigl(
        \l|R_2|^2
        +2\l\Re(\psi^\e\xwl{R_2})
        +\l|\psi^\e|^2
        \bigr)
        =0, \\
        \xwt{R_1}(0)=0.
    \end{dcases}
\end{align*}

% By \eqref{init-R1R2}, we have

\subsection{Local existence and continuation} %\label{sec:local-existence}
We first establish local existence and uniqueness for the remainder system \eqref{eq-for-R1R2}, together with a continuation criterion.
The proof uses the profile bounds from Section~2, the linear estimates, and Lemma~\ref{lem:2}.

\begin{proposition} \label{prop:local-original}
    Under the assumptions of Theorem~\ref{thm:1}, the system \eqref{eq-for-R1R2} with initial data \eqref{init-R1R2} admits a unique maximal solution in the class described below, defined on $[0,T^\e)$ for some $T^\e \in (0,\infty]$. On every compact interval $I\subset[0,T^\e)$, the solution satisfies
    \begin{align}\label{local-class}
        \begin{split}
             & \xwt{R_1}\in C(I;H^{\s+\frac12}(\R^d)),
            \qquad R_2\in C(I;H^\s(\R^d)),                     \\
             & |\nabla|^\g\xwt{R_1},\ |\nabla|^\g R_2\in X(I),
            \qquad \g\in\{0,\s\}.
        \end{split}
    \end{align}
    If $T^\e<\infty$, then
    \begin{align}\label{continuation-criterion}
        \sum_{\g\in\{0,\s\}}
        \left(
        \||\nabla|^\g\xwt{R_1}\|_{X([0,T^\e))}
        +\||\nabla|^\g R_2\|_{X([0,T^\e))}
        \right)=\infty.
    \end{align}
\end{proposition}

\begin{proof}
    We first construct a local solution from prescribed data at an arbitrary time $t_0\ge0$, with
    \begin{align*}
        \xwt{R_1}(t_0)\in H^{\s+\frac12}(\R^d),
        \qquad R_2(t_0)\in H^\s(\R^d).
    \end{align*}
    For functions $u$ and $v$, define
    \begin{align*}
        F_1(v)   & :=\l|v+\psi^\e|^2,                                      \\
        F_2(u,v) & :=-\l\left(2\Re(e^{it}\phi^\e)+\Im u\right)(v+\psi^\e).
    \end{align*}
    These are precisely the forcing terms in the equations for $\xwt{R_1}$ and $R_2$ when $(u,v)=(\xwt{R_1},R_2)$ in \eqref{eq-for-R1R2}. On $I=[t_0,t_0+h]$, define the Duhamel map $\mathcal T=(\mathcal T_1,\mathcal T_2)$ by
    \begin{align*}
        \mathcal T_1(u,v)(t)
         & :=e^{i(t-t_0)\lna}\xwt{R_1}(t_0)
        +\int_{t_0}^t e^{i(t-s)\lna}\lna^{-1}F_1(v)(s)\,\D{s}, \\
        \mathcal T_2(u,v)(t)
         & :=e^{i(t-t_0)\Delta}R_2(t_0)
        -i\int_{t_0}^t e^{i(t-s)\Delta}F_2(u,v)(s)\,\D{s}.
    \end{align*}
    Choose $M>0$ such that
    \begin{align*}
        M>2\sum_{\g\in\{0,\s\}}
        \left(
        \|\lna^{\frac12}|\nabla|^\g\xwt{R_1}(t_0)\|_{L^2}
        +\||\nabla|^\g R_2(t_0)\|_{L^2}
        \right).
    \end{align*}
    By unitarity, for every $t\in  I =[t_{0}, t_{0} + h]$,
    \begin{align*}
         & \sum_{\g\in\{0,\s\}}
        \Bigl(
        \|\lna^{\frac12}|\nabla|^\g e^{i(t-t_0)\lna}\xwt{R_1}(t_0)\|_{L^2}
        +\||\nabla|^\g e^{i(t-t_0)\Delta}R_2(t_0)\|_{L^2}
        \Bigr)
        <\frac{M}{2}.
    \end{align*}

    Let $\rho>0$ be sufficiently small, to be chosen below.
    By the linear estimates \eqref{est:s-3} and \eqref{est:kg-3}, together with the profile estimates \eqref{est:phi-0} and \eqref{est:psi-0}, we can choose $h>0$ sufficiently small that
    \begin{align*}
        \sum_{\g\in\{0,\s\}}\Bigl(
         & \||\nabla|^\g e^{i(t-t_0)\lna}\xwt{R_1}(t_0)\|_{X(I)}
        +\||\nabla|^\g e^{i(t-t_0)\Delta}R_2(t_0)\|_{X(I)}       \\
         & +\||\nabla|^\g\phi^\e\|_{X(I)}
        +\||\nabla|^\g\psi^\e\|_{X(I)}
        \Bigr)\le\frac{\rho}{4}.
    \end{align*}
    Consider the set
    \begin{align*}
        \mathcal B = \Big\{ (u,v)\in C(I;H^{\s+\frac12}(\R^d)\times H^\s(\R^d)): & \sum_{\g\in\{0,\s\}}
        \big(
        \|\lna^{\frac12}|\nabla|^\g u\|_{L_t^\infty L_x^2(I)}
        +\||\nabla|^\g v\|_{L_t^\infty L_x^2(I)}
        \big) \le M,                                                                                    \\
                                                                                 & \sum_{\g\in\{0,\s\}}
        \big(
        \||\nabla|^\g u\|_{X(I)}
        +\||\nabla|^\g v\|_{X(I)}
        \big)\le\rho. \Big\}
    \end{align*}

    Applying   Lemma~\ref{lem:2}  to $F_1(v)=\l(v+\psi^\e)\xwl{(v+\psi^\e)}$ gives
    \begin{align*}
        \||\nabla|^\g F_1(v)\|_{L_{tx}^{\frac{2(d+2)}{d+4}}(I)}
        \le C
        \||\nabla|^\g(v+\psi^\e)\|_{X(I)}
        \||\nabla|^{s_c}(v+\psi^\e)\|_{X(I)}, \quad \mbox{for $\g\in\{0,\s\}$.}
    \end{align*}
    Since $\s \geq s_{c} = \frac{d}{2}-2$, interpolation yields
    \begin{align*}
        \||\nabla|^{s_c}(v+\psi^\e)\|_{X(I)}
        \le C
        \sum_{\g\in\{0,\s\}}
        \||\nabla|^\g(v+\psi^\e)\|_{X(I)}.
    \end{align*}
    Consequently,
    \begin{align*}
        \sum_{\g\in\{0,\s\}}
        \||\nabla|^\g F_1(v)\|_{L_{tx}^{\frac{2(d+2)}{d+4}}(I)}
        \le C
        \big(
        \sum_{\g\in\{0,\s\} }
        \||\nabla|^\g(v+\psi^\e)\|_{X(I)}
        \big)^2.
    \end{align*}
    By the definition of $\mathcal B$ and the choice of $h$,
    \begin{align*}
        \sum_{\g\in\{0,\s\}}
        \||\nabla|^\g(v+\psi^\e)\|_{X(I)}
        \le \rho+\frac{\rho}{4},
    \end{align*}
    which gives
    \begin{align*}
        \sum_{\g\in\{0,\s\}}
        \||\nabla|^\g F_1(v)\|_{L_{tx}^{\frac{2(d+2)}{d+4}}(I)}
        \le C\rho^2.
    \end{align*}
    The same argument applies to $F_2(u,v)$:
    \begin{align*}
        \sum_{\g\in\{0,\s\}}
        \||\nabla|^\g F_2(u,v)\|_{L_{tx}^{\frac{2(d+2)}{d+4}}(I)}
        \le C\rho^2.
    \end{align*}

    Choose $\rho>0$ so that
    \begin{align*} %\label{def:rho}
        C\rho\le\frac14,
        \qquad C\rho^2\le\frac{M}{2}.
    \end{align*}
    Thus,  by \eqref{est:s-3}, \eqref{est:kg-4}, and the preceding bounds for $F_1$ and $F_2$, we obtain
    \begin{align*}
         & \sum_{\g\in\{0,\s\}}
        \Bigl(
        \|\lna^{\frac12}|\nabla|^\g\mathcal T_1(u,v)\|_{L_t^\infty L_x^2(I)}
        +\||\nabla|^\g\mathcal T_2(u,v)\|_{L_t^\infty L_x^2(I)}
        \Bigr)        \le\frac{M}{2}+C\rho^2\le M.
    \end{align*}
    Similarly,
    \begin{align*}
         & \sum_{\g\in\{0,\s\}}
        \Bigl(
        \||\nabla|^\g\mathcal T_1(u,v)\|_{X(I)}
        +\||\nabla|^\g\mathcal T_2(u,v)\|_{X(I)}
        \Bigr)                   \le\frac{\rho}{4}+C\rho^2
        \le\frac{\rho}{2}.
    \end{align*}
    The same linear estimates also give
    \begin{align*}
        \mathcal T_1(u,v)\in C(I;H^{\s+\frac12}(\R^d)),
        \qquad
        \mathcal T_2(u,v)\in C(I;H^\s(\R^d)).
    \end{align*}
    Thus $\mathcal T$ maps $\mathcal B$ into itself.

    We equip $\mathcal B$ with the norm
    \begin{align*}
        \|  (u, v)\|_{\mathcal B} : =        \sum_{\g\in\{0,\s\}}
        \bigl(
        \| |\nabla|^\g u\|_{Y^{\frac12}(I)}   +\||\nabla|^\g v\|_{Y^0(I)}
        \bigr).
    \end{align*}
    For $(u_1,v_1),(u_2,v_2)\in \mathcal B$, applying Lemma~\ref{lem:2} and the fact $\s \geq s_{c}$ gives
    \begin{align*}
         & \| \mathcal T(u_1,v_1) - \mathcal T(u_2,v_2) \|_{{\mathcal B}} \\
         & \le  C \sum_{\g\in\{0,\s\}}
        \Bigl(
        \||\nabla|^\g(F_1(v_1)-F_1(v_2))\|_{L_{tx}^{\frac{2(d+2)}{d+4}}(I)}
        +\||\nabla|^\g(F_2(u_1,v_1)-F_2(u_2,v_2))\|_{L_{tx}^{\frac{2(d+2)}{d+4}}(I)}
        \Bigr).
    \end{align*}

    We again apply Lemma~\ref{lem:2} to
    \begin{align*}
        F_1(v_1)-F_1(v_2)
        = & \l(v_1-v_2)\xwl{(v_1+\psi^\e)}
        +\l(v_2+\psi^\e)\xwl{(v_1-v_2)},
    \end{align*}
    and obtain
    \begin{align*}
         & \sum_{\g\in\{0,\s\}}
        \||\nabla|^\g(F_1(v_1)-F_1(v_2))\|_{L_{tx}^{\frac{2(d+2)}{d+4}}(I)} \\
         & \quad\le C
        \big(
        \sum_{j=1}^2\sum_{\g\in\{0,\s\}}
        \||\nabla|^\g(v_j+\psi^\e)\|_{X(I)}
        \big)
        \big(
        \sum_{\g\in\{0,\s\}}
        \||\nabla|^\g(v_1-v_2)\|_{X(I)}
        \big)                                                               \\
         & \quad\le C\rho
        \sum_{\g\in\{0,\s\}}
        \||\nabla|^\g(v_1-v_2)\|_{X(I)},
    \end{align*}
    where we used the definition of $\mathcal B$ and the choice of $h$ to bound
    \begin{align*}
        \sum_{\g\in\{0,\s\}}
        \||\nabla|^\g(v_j+\psi^\e)\|_{X(I)}
        \le\rho+\frac{\rho}{4} ,
        \qquad j=1,2.
    \end{align*}
    The same estimate holds for $F_2(u_1,v_1)-F_2(u_2,v_2)$. Consequently,
    \begin{align*}
         & \bigl\| \mathcal T(u_1,v_1) - \mathcal T(u_2,v_2)\bigr\|_{\mathcal B} \le C\rho\,
        \bigl\|  (u_1,v_1) - (u_2,v_2) \bigr\|_{\mathcal B} \le \frac14
        \bigl\| (u_1,v_1) - (u_2,v_2) \bigr\|_{\mathcal B} .
    \end{align*}
    Therefore $\mathcal T$ is a contraction on $\mathcal B$, and the contraction mapping theorem yields a unique fixed point $(u,v)$ which is the solution to \eqref{eq-for-R1R2}--\eqref{init-R1R2}.

    It remains to prove the continuation criterion. Suppose that $T^\e<\infty$ and that the sum in \eqref{continuation-criterion} is finite. Write $(u,v)=(\xwt{R_1},R_2)$. Interpolation, the profile bounds, and Lemma~\ref{lem:2} imply
    \begin{align*}
        |\nabla|^\g F_1(v),\
        |\nabla|^\g F_2(u,v)
        \in L_{tx}^{\frac{2(d+2)}{d+4}}([0,T^\e)),
        \qquad \g\in\{0,\s\}.
    \end{align*}
    For $0\le t<t'<T^\e$, the Duhamel formulas and the linear estimates give
    \begin{align*}
         & \|e^{-it'\lna}u(t')-e^{-it\lna}u(t)\|_{H^{\s+\frac12}}
        +\|e^{-it'\Delta}v(t')-e^{-it\Delta}v(t)\|_{H^\s}         \\
         & \quad\le C\sum_{\g\in\{0,\s\}}
        \big(
        \||\nabla|^\g F_1(v)\|_{L_{tx}^{\frac{2(d+2)}{d+4}}([t,t'])}
        +\||\nabla|^\g F_2(u,v)\|_{L_{tx}^{\frac{2(d+2)}{d+4}}([t,t'])}
        \big).
    \end{align*}
    The right-hand side tends to zero as $t,t'\to T^\e$. Hence
    \begin{align*}
        e^{-it\lna}u(t),
        \qquad
        e^{-it\Delta}v(t)
    \end{align*}
    are Cauchy as $t\to T^\e$ in
    $H^{\s+\frac12}(\R^d)$ and $H^\s(\R^d)$, respectively.

    By completeness, these two families converge. The strong continuity of the group then gives
    \begin{align*}
        \lim_{t \to T^\e}u(t)
         & =
        e^{i T^\e\lna}
        \lim_{t\to T^\e}e^{-it\lna}u(t)
        \quad\text{in }H^{\s+\frac12}(\R^d), \\
        \lim_{t\to T^\e}v(t)
         & =
        e^{i T^\e\Delta}
        \lim_{t\to T^\e}e^{-it\Delta}v(t)
        \quad\text{in }H^\s(\R^d).
    \end{align*}
    We denote these  limits by $u(T^\e)$ and $v(T^\e)$, respectively.
    Using them as initial data at $t_0=T^\e$, the local construction extends the solution beyond $T^\e$. This contradicts maximality and proves \eqref{continuation-criterion}.

\end{proof}

\subsection{A finite time partition and smallness of the profiles}

The a priori estimates below require small $X$ norms of $|\nabla|^{s_c}\phi^\e$ and $|\nabla|^{s_c}\psi^\e$ on each time interval. The global profile bounds in Section 2 allow us to construct a finite partition of $\R_+$ with a number of intervals bounded independently of $\e$.

\begin{proposition} \label{prop:time-partition}
    Under the assumptions of Theorem~\ref{thm:1}, for every $\de>0$ and $0<\e\le1$, there exists a finite partition
    \begin{align*} %\label{ti}
        0 = t_0<t_1<\cdots < t_n=\infty,
    \end{align*}
    such that, with $I_k=[t_k,t_{k+1})$,
    \begin{align}\label{1-4zw}
        \||\nabla|^{s_c}\phi^\e\|_{X(I_k)}\le\de,\qquad
        \||\nabla|^{s_c}\psi^\e\|_{X(I_k)}\le\de,
        \qquad k=0,\ldots,n-1.
    \end{align}
    Moreover,
    \begin{align}\label{partition-number}
        n \le \lfloor N_\de \rfloor + 1, \qquad  N_\de:=1+2\sum_{j=1}^2
        \left(\frac{8C C_0}{\de}\right)^{p_j}.
    \end{align}
\end{proposition}

\begin{proof}
    Taking $\g=s_c$ in \eqref{est:phi-0} and \eqref{est:psi-0}, we obtain
    \begin{align*}
        \||\nabla|^{s_c}\phi^\e\|_{X(\R_+)}\le C C_0,\qquad
        \||\nabla|^{s_c}\psi^\e\|_{X(\R_+)}\le C C_0.
    \end{align*}
    Set $t_0=0$. Suppose that $t_k$ has been chosen. If
    \begin{align*}
        \||\nabla|^{s_c}\phi^\e\|_{X([t_k,\infty))}
        +\||\nabla|^{s_c}\psi^\e\|_{X([t_k,\infty))}
        \le\frac{\de}{2},
    \end{align*}
    take $[t_k,\infty)$ as the final interval. Otherwise, by the definition of $X$ in \eqref{def:norm-X},  we can choose a finite $t_{k+1}>t_k$ such that
    \ba\label{cond:de/2}
    \||\nabla|^{s_c}\phi^\e\|_{X(I_k)}
    +\||\nabla|^{s_c}\psi^\e\|_{X(I_k)}
    =\frac{\de}{2}.
    \ea
    By \eqref{def:norm-X}, the left-hand side is the sum of four component norms. Thus at least one of them is at least $\frac{\de}{8}$, and hence
    \begin{align*}
        1\le\sum_{j=1}^2\left(\frac8\de\right)^{p_j}
        \big(
        \||\nabla|^{s_c}\phi^\e\|_{L_t^{p_j}L_x^{q_j}(I_k)}^{p_j}
        +\||\nabla|^{s_c}\psi^\e\|_{L_t^{p_j}L_x^{q_j}(I_k)}^{p_j}
        \big).
    \end{align*}
    If $m$ finite intervals have been constructed in this way, summing over these disjoint intervals gives
    \ba\label{m-interval-number}
    m & \le\sum_{j=1}^2\left(\frac8\de\right)^{p_j}
    \sum_{k=0}^{m-1}
    \big(
    \||\nabla|^{s_c}\phi^\e\|_{L_t^{p_j}L_x^{q_j}(I_k)}^{p_j}
    +\||\nabla|^{s_c}\psi^\e\|_{L_t^{p_j}L_x^{q_j}(I_k)}^{p_j}
    \big)                                                       \\
    & \le\sum_{j=1}^2\left(\frac8\de\right)^{p_j}
    \big(
    \||\nabla|^{s_c}\phi^\e\|_{L_t^{p_j}L_x^{q_j}(\R_+)}^{p_j}
    +\||\nabla|^{s_c}\psi^\e\|_{L_t^{p_j}L_x^{q_j}(\R_+)}^{p_j}
    \big)                                                       \\
    & \le 2\sum_{j=1}^2\left(\frac{8 C C_0}{\de}\right)^{p_j}.
    \ea
    Therefore, the construction must stop after finitely many steps with a final interval $[t_{n-1},\infty)$. The bound on the number of intervals in \eqref{partition-number} follows immediately from \eqref{m-interval-number}. The estimates in \eqref{1-4zw} follow from the choice of $I_{k}$, see \eqref{cond:de/2}.
\end{proof}

\subsection{Estimate for the wave-type remainder}

We now derive a priori estimates on $I_k = [t_{k}, t_{k+1})$ defined in Section 3.2, assuming that the solution exists on this interval. The same estimates will hold with $I_k$ replaced by $[t_k,\tau]$ for every
\begin{align*}
    t_k<\tau<\min\{t_{k+1},T^\e\},
\end{align*}
with estimate constants independent of $\tau$.

For $t\in I_k$, Duhamel's formula gives
\begin{align} \label{duhamel-R1}
    \begin{split}
        \xwt{R_1}(t) & = e^{i(t-t_k) \lna} \xwt{R_{1}} (t_k)   +  \int_{t_k}^t e^{i(t-s) \lna} \lna^{-1} \big( \l |R_2|^2 + 2 \l \Re( \psi^\e \xwl{R_2}) + \l | \psi^\e|^2 \big)(s) \D{s}.
    \end{split}
\end{align}

To estimate the forcing $|\psi^\e|^2$, we integrate the corresponding term in \eqref{duhamel-R1} by parts in time, using
$
    \pa_s e^{i(t-s)\lna}=-i\lna e^{i(t-s)\lna}
$, and derive that
\begin{align*}
    \begin{split}
        \xwt{R_1}(t)
        = & e^{i(t-t_k)\lna}\xwt{R_1}(t_k)
        +i\l\lna^{-2}|\psi^\e|^2(t)
        -i\l e^{i(t-t_k)\lna}\lna^{-2}|\psi^\e|^2(t_k)     \\
          & +\l\int_{t_k}^t e^{i(t-s)\lna}\lna^{-1}
        \bigl(|R_2|^2+2\Re(\psi^\e\xwl{R_2})\bigr)(s)\D{s} \\
          & -i\l\int_{t_k}^t e^{i(t-s)\lna}\lna^{-2}
        2\Re\bigl((\pa_s\psi^\e)\xwl{\psi^\e}\bigr) \D{s},
    \end{split}
\end{align*}
where we used the fact $
    \pa_t |\psi^\e|^2= 2\Re\bigl((\pa_t\psi^\e)\xwl{\psi^\e}\bigr).
$

To estimate the terms evaluated at the time interval endpoints $t_k$ and $t$, choose
% To estimate the terms at ending points $t$ and $t_{k}$, choose
\be\label{choice-eta}
0<\eta<\min \big\{2,\frac d2,\frac{d}{q_1},\frac{d}{q_2}\big\}.
\ee
Since $\s\ge s_c=\frac d2-2$, we have
\begin{align*}
    \g+\eta\le\s+2,\qquad \frac d2-\eta\le\s+2,
\end{align*}
for every $0\le\g\le\s$. Applying Lemma~\ref{ineq:eta} with spatial exponents $2,q_1,q_2$,  together with  H\"older's inequality in time and \eqref{est:psi-0}, gives
\begin{align}\label{est:phil}
    \begin{split}
        \||\nabla|^\g|\psi^\e|^2\|_{Y^0(\R_+)}
         & \le C\||\nabla|^{\g+\eta}\psi^\e\|_{Y^0(\R_+)}
        \||\nabla|^{\frac d2-\eta}\psi^\e\|_{L_t^\infty L_x^2(\R_+)} \\
         & \le C C_0^2 \e^{\g+\eta-s_c}\e^{\frac d2-\eta-s_c}
        =C C_0^2 \e^{\g+2-s_c}.
    \end{split}
\end{align}

By the boundedness of $\lna^{-\frac{3}{2}}$ on $L_x^2(\R^d)$ and $\lna^{-2}$ on the spatial $L^r$ components of $X$, and by \eqref{est:kg-3}, we have
\begin{align*}
    \bigl\||\nabla|^\g\lna^{-2}|\psi^\e|^2\bigr\|_{Y^{\frac{1}{2}}(I_k)}                       & \le C\bigl\||\nabla|^\g|\psi^\e|^2\bigr\|_{Y^0(I_k)},                                                       \\
    \bigl\||\nabla|^\g e^{i(t-t_k)\lna}\lna^{-2}|\psi^\e(t_k)|^2\bigr\|_{Y^{\frac{1}{2}}(I_k)} & \le C\bigl\||\nabla|^\g|\psi^\e(t_k)|^2\bigr\|_{L_x^2}\le C\bigl\||\nabla|^\g|\psi^\e|^2\bigr\|_{Y^0(I_k)}.
\end{align*}

Taking the $Y^{\frac{1}{2}}$ norm on $I_k$ and applying \eqref{est:kg-3} and \eqref{est:kg-4}, we obtain
\begin{align} \label{est:R1-g}
    \begin{split}
        \| |\nabla|^\g \xwt{R_1} \|_{Y^{\frac{1}{2}}(I_k)} \le & C \| \lna^{\frac{1}{2}} |\nabla|^\g \xwt{R_1} (t_k) \|_{L_x^2(\R^d)} + C   \| |\nabla|^\g |\psi^\e |^2 \|_{Y^0(I_k) }  + C   \| |\nabla|^\g  |R_2|^2 \|_{L_{t x}^{\frac{2(d+2)}{d+4}}(I_k) }     \\
                                                               & \quad + C   \| |\nabla|^\g \Re( \psi^\e \xwl{R_2})   \|_{L_{t x}^{\frac{2(d+2)}{d+4}}(I_k) }  + C \| |\nabla|^\g \big((\pa_t \psi^\e )\xwl{\psi^\e}\big) \|_{L_{tx}^{\frac{2(d+2)}{d+4}}(I_k) }.
    \end{split}
\end{align}
Lemma~\ref{lem:2} controls the terms $|R_2|^2$ and $\psi^\e\overline{R_2}$ by the corresponding products of $X(I_k)$ norms. Applying the same lemma together with \eqref{est:psi-0} and \eqref{est:psi-1}, we also obtain
\begin{align*}
    \||\nabla|^\g((\pa_t\psi^\e)\overline{\psi^\e})\|_{L_{tx}^{\frac{2(d+2)}{d+4}}(I_k)}
     & \le C\||\nabla|^\g\pa_t\psi^\e\|_{X(I_k)}\||\nabla|^{s_c}\psi^\e\|_{X(I_k)}
    +C\||\nabla|^{s_c}\pa_t\psi^\e\|_{X(I_k)}\||\nabla|^\g\psi^\e\|_{X(I_k)}       \\
     & \le C C_0^2 \e^{\g+2-s_c}.
\end{align*}
Combining this with \eqref{est:R1-g} and \eqref{est:phil}, we have
\begin{align} \label{R1-general}
    \begin{split}
        \| |\nabla|^\g \xwt{R_1} \|_{ Y^{\frac{1}{2}} (I_k)} & \leq  C \| \lna^{\frac{1}{2} } |\nabla|^{\g} \xwt{R_1} (t_k) \|_{L^2 (\R^d)} + C \| |\nabla|^\g R_2 \|_{X (I_k)} \| |\nabla|^{s_c} R_2 \|_{X (I_k)}                \\
                                                             & \quad  + C \| |\nabla|^\g \psi^\e \|_{ X (I_k) } \| |\nabla|^{s_c} R_2 \|_{ X(I_k) } + C \| |\nabla|^{s_c} \psi^\e \|_{X (I_k) } \| |\nabla|^\g R_2 \|_{ X (I_k) } \\
                                                             & \quad  + C C_0^2 \e^{\g+2-s_c}.
    \end{split}
\end{align}
Taking $\g=s_c$ in \eqref{R1-general} and using \eqref{1-4zw}, we obtain
\begin{align} \label{est-R1}
    \begin{split}
        \| |\nabla|^{s_c} \xwt{R_1} \|_{ Y^{\frac{1}{2}}(I_k) }
         & \leq  C \|\lna^{\frac{1}{2}} |\nabla|^{s_c}  \xwt{R_1} (t_k) \|_{L^2 (\R^d)} + C \| |\nabla|^{s_c} R_2\|_{X(I_k)}^2 \\
         & \qquad +  C \de \| |\nabla|^{s_c} R_2 \|_{X(I_k)}   + C C_0^2 \e^2.
    \end{split}
\end{align}

\subsection{Estimate for the Schr\"odinger-type remainder} %\label{sec-R2}

We next estimate the Schr\"odinger-type remainder $R_2$ on $I_k=[t_k,t_{k+1})$ defined in Section 3.2. The same estimates hold on $[t_k,\tau]$ with $t_k<\tau<\min\{t_{k+1},T^\e\}$.
The terms containing $R_1$ or $R_2$ are estimated by Lemma~\ref{lem:2}.
For the forcing term $2\l\Re(e^{it}\phi^\e)\psi^\e$, we separate the low and high spatial frequencies. We treat the low frequencies by integration by parts and estimate the high frequencies using two additional spatial derivatives.

For $t\in I_k$, Duhamel's formula gives
\begin{align} \label{duhamel-R2}
    \begin{split}
        R_2(t) & = e^{i(t-t_k) \Delta } R_2(t_k)                                                                                                                           +  2 i \l  \int_{t_k}^t  e^{i(t-s)\Delta} \big ( P_{\leq \frac{1}{2}} (\Re(e^{i s}  \phi^\e )  \psi^\e ) \big) \D{s} \\
               & \quad + i \int_{t_k}^t e^{i(t-s)\Delta } \big( 2 \l \Re ( e^{i s}  \phi^\e ) R_2 + \l R_1 \psi^\e + \l R_1 R_2 + 2\l  P_{> \frac{1}{2}} (\Re(e^{i s}  \phi^\e )  \psi^\e ) \big) \D{s}.
    \end{split}
\end{align}
Here $R_2(0)=0$ by \eqref{init-R1R2}.

The low-frequency forcing in \eqref{duhamel-R2} can be written as
\begin{align*}
     & 2i\l\int_{t_k}^t e^{i(t-s)\Delta}P_{\le\frac{1}{2}}
    \bigl(\Re(e^{is}\phi^\e)\psi^\e\bigr)\D{s} = i\l e^{it\Delta}\int_{t_k}^t
    \bigl(e^{-is(\Delta-1)}P_{\le\frac{1}{2}}(\phi^\e\psi^\e)
    +e^{-is(\Delta+1)}P_{\le\frac{1}{2}}(\xwl{\phi^\e}\psi^\e)\bigr) \D{s}.
\end{align*}
The cutoff is supported in $|\xi|\le\frac{1}{2}$, so both $-1-|\xi|^2$ and $1-|\xi|^2$ are bounded away from zero there. In particular, $(\Delta\pm 1)^{-1}P_{\le\frac{1}{2}}$ are bounded multipliers on $L_x^r$ for $1<r<\infty$.

Integration by parts gives
\begin{align}\label{ibp-R2}
    \begin{split}
         & i\l e^{it\Delta}\int_{t_k}^t
        e^{-is(\Delta-1)}P_{\le\frac{1}{2}}(\phi^\e\psi^\e)(s)\D{s}                                               \\
         & \quad = -\l e^{it}(\Delta-1)^{-1}P_{\le\frac{1}{2}}(\phi^\e\psi^\e)(t)  +\l e^{i(t-t_k)\Delta}e^{it_k}
        (\Delta-1)^{-1}P_{\le\frac{1}{2}}(\phi^\e\psi^\e)(t_k)                                                    \\
         & \qquad+\l\int_{t_k}^t e^{i(t-s)\Delta}e^{is}
        (\Delta-1)^{-1}P_{\le\frac{1}{2}}\pa_s(\phi^\e\psi^\e)(s) \D{s},
    \end{split}
\end{align}
and
\begin{align}\label{ibp-R2-new}
    \begin{split}
         & i\l e^{it\Delta}\int_{t_k}^t
        e^{-is(\Delta+1)}P_{\le\frac{1}{2}}(\overline {\phi^\e}\psi^\e)(s)\D{s}                                                  \\
         & \quad = -\l e^{-it}(\Delta+1)^{-1}P_{\le\frac{1}{2}}(\overline {\phi^\e} \psi^\e)(t)  +\l e^{i(t-t_k)\Delta}e^{-it_k}
        (\Delta+1)^{-1}P_{\le\frac{1}{2}}( \overline {\phi^\e} \psi^\e)(t_k)                                                     \\
         & \qquad+\l\int_{t_k}^t e^{i(t-s)\Delta}e^{-i s}
        (\Delta+1)^{-1}P_{\le\frac{1}{2}}\pa_s( \overline {\phi^\e} \psi^\e)(s)\D{s}.
    \end{split}
\end{align}

The boundedness of $(\Delta\mp 1)^{-1} P_{\le \frac12}$, together with \eqref{est:s-3}, \eqref{ibp-R2}--\eqref{ibp-R2-new}, gives the following estimate: for $0\le\g\le\s$,
\begin{align}\label{R2-before-high}
    \begin{split}
         & \||\nabla|^\g R_2\|_{Y^0(I_k)}    \le  C \||\nabla|^\g R_2(t_k)\|_{L^2}
        + C \||\nabla|^\g  (\phi^\e\psi^\e, \xwl{\phi^\e}\psi^\e)\|_{Y^0(I_k)}     \\
         & \qquad   + C
        \||\nabla|^\g\pa_t   (\phi^\e\psi^\e, \xwl{\phi^\e}\psi^\e) \|_{L_{tx}^{\frac{2(d+2)}{d+4}}(I_k)}  +  C \||\nabla|^\g P_{>\frac{1}{2}} \bigl(\Re(e^{it}\phi^\e)\psi^\e\bigr)
        \|_{L_{tx}^{\frac{2(d+2)}{d+4}}(I_k)}                                      \\
         & \qquad + C \||\nabla|^\g\bigl(\Re(e^{it}\phi^\e)R_2
        +R_1\psi^\e+R_1 R_2\bigr)\|_{L_{tx}^{\frac{2(d+2)}{d+4}}(I_k)}   .
    \end{split}
\end{align}

Choosing $\eta$ as in \eqref{choice-eta} and applying Lemma~\ref{ineq:eta} together with the profile bounds \eqref{est:phi-0} and \eqref{est:psi-0}, we obtain, for $0\le \g \le \s$,
\begin{align*} %\label{est:psil}
    \begin{split}
        \||\nabla|^\g(\phi^\e\psi^\e)\|_{Y^0(\R_+)} +\||\nabla|^\g(\overline{\phi^\e}\psi^\e)\|_{Y^0(\R_+)} & \le C\||\nabla|^{\g+\eta}\phi^\e\|_{Y^0(\R_+)}
        \||\nabla|^{\frac d2-\eta}\psi^\e\|_{L_t^\infty L_x^2(\R_+)}                                                                                                               \\
                                                                                                            & \qquad+C\||\nabla|^{\frac d2-\eta}\phi^\e\|_{L_t^\infty L_x^2(\R_+)}
        \||\nabla|^{\g+\eta}\psi^\e\|_{Y^0(\R_+)}                                                                                                                                  \\
                                                                                                            & \le C C_0^2 \e^{\g+2-s_c}.
    \end{split}
\end{align*}

Write $\pa_t(\phi^\e\psi^\e)=(\pa_t\phi^\e)\psi^\e+\phi^\e\pa_t\psi^\e$. By Lemma~\ref{lem:2}, \eqref{est:phi-1}, and \eqref{est:psi-0}, we have
\begin{align*}
    \||\nabla|^\g((\pa_t\phi^\e)\psi^\e)\|_{L_{tx}^{\frac{2(d+2)}{d+4}}(I_k)}
     & \le C\||\nabla|^\g\pa_t\phi^\e\|_{X(I_k)}\||\nabla|^{s_c}\psi^\e\|_{X(I_k)}
    +C\||\nabla|^{s_c}\pa_t\phi^\e\|_{X(I_k)}\||\nabla|^\g\psi^\e\|_{X(I_k)}       \\
     & \le C C_0^2 \e^{\g+2-s_c}.
\end{align*}
The term $\phi^\e\pa_t\psi^\e$ has the same bound by \eqref{est:phi-0} and \eqref{est:psi-1}. The condition $s_c \le \s$ ensures that both time-derivative estimates are available. Complex conjugation does not change these estimates. Thus
\begin{align}\label{R2-2}
    \||\nabla|^\g\pa_t (\phi^\e\psi^\e, \xwl{\phi^\e}\psi^\e) \|_{L_{tx}^{\frac{2(d+2)}{d+4}}(I_k)}
    \le C C_0^2 \e^{\g+2-s_c},\qquad 0\le\g\le\s.
\end{align}
It remains to estimate the high-frequency term in \eqref{R2-before-high}. Since $P_{>\frac{1}{2}}$ vanishes on $|\xi|\le\frac{1}{4}$, $P_{>\frac{1}{2}}|\nabla|^{-2}$ is bounded on $L^r(\R^d)$, $1<r<\infty$.
Applying Lemma~\ref{lem:2} at order $\g+2$ and using \eqref{est:phi-0} and \eqref{est:psi-0}, we obtain, for $0\le\g\le\s$,
\begin{align*}
     & \||\nabla|^\g P_{>\frac{1}{2}}  (\Re(e^{it} \phi^\e )\psi^\e)
    \|_{L_{tx}^{\frac{2(d+2)}{d+4}}(I_k)}           \le C \||\nabla|^{\g+2}(\phi^\e\psi^\e)
    \|_{L_{tx}^{\frac{2(d+2)}{d+4}}(I_k)}                            \\
     & \quad\le C \||\nabla|^{\g+2}\phi^\e\|_{X(I_k)}
    \||\nabla|^{s_c}\psi^\e\|_{X(I_k)}
    + C \||\nabla|^{s_c}\phi^\e\|_{X(I_k)}
    \||\nabla|^{\g+2}\psi^\e\|_{X(I_k)}                              \\
     & \quad\le C C_0^2 \e^{\g+2-s_c}.
\end{align*}
At $\g=\s$, we need initial regularity $\phi_0,\phi_1\in H^{\s+\frac52}(\R^d)$ and $\psi_0\in H^{\s+2}(\R^d)$, as assumed in \eqref{init:0}. The same calculation applies to $\xwl{\phi^\e}\psi^\e$, giving
\begin{align}\label{R2-high}
    \||\nabla|^\g P_{>\frac{1}{2}}\bigl(\Re(e^{it}\phi^\e)\psi^\e\bigr)
    \|_{L_{tx}^{\frac{2(d+2)}{d+4}}(I_k)}
    \le C C_0^2 \e^{\g+2-s_c}.
\end{align}

Combining \eqref{R2-before-high}, \eqref{R2-2}, \eqref{R2-high}, and the quadratic remainder bounds from Lemma~\ref{lem:2}, we obtain
\begin{align} \label{R2-general}
    \begin{split}
             & \| |\nabla|^\g R_2 \|_{Y^0 (I_k) }                                                                                                                                                                          \\
        \leq & C \| |\nabla|^{\g}  R_2(t_k) \|_{L_x^2(\R^d)} +  C \| |\nabla|^{s_c}   \phi^\e \|_{ X(I_k) } \| |\nabla|^\g R_2 \|_{ X(I_k) } + C \| |\nabla|^\g \phi^\e \|_{ X(I_k) } \| |\nabla|^{s_c}  R_2 \|_{ X(I_k) } \\
             & + C \| |\nabla|^{s_c} \psi^\e \|_{ X(I_k) } \| |\nabla|^\g R_1 \|_{ X(I_k) } + C \| |\nabla|^\g \psi^\e \|_{ X(I_k) }   \| |\nabla|^{s_c}  R_1 \|_{ X(I_k) }                                                \\
             & + C  \| |\nabla|^\g R_1\|_{ X(I_k) } \| |\nabla|^{s_c} R_2 \|_{X(I_k) } + C \| |\nabla|^{s_c} R_1 \|_{ X(I_k) } \| |\nabla|^\g R_2 \|_{ X(I_k) }  +  C C_0^2 \e^{\g+2-s_c}.
    \end{split}
\end{align}
Setting $\g = s_c$ in \eqref{R2-general} gives
\begin{align} \label{est-R2}
    \begin{split}
        \| |\nabla|^{s_c} R_2 \|_{Y^0(I_k) } \leq & C \| |\nabla|^{s_c}  R_2(t_k) \|_{L_x^2(\R^d)} + 2 C \de \| |\nabla|^{s_c} R_2 \|_{X(I_k) }   +  2 C \de  \| |\nabla|^{s_c} R_1 \|_{ X(I_k)} \\
                                                  & + 2 C \| |\nabla|^{s_c} R_1 \|_{ X(I_k)} \| |\nabla|^{s_c} R_2 \|_{X(I_k)}  +C C_0^2 \e^2.
    \end{split}
\end{align}

\subsection{Global error estimates} %\label{sec:global-estimates}

In this subsection, under the assumptions of Theorem~\ref{thm:1}, we will establish uniform error bounds for $(\xwt{R_1},R_2)$ on each interval $I_k$ and then prove global existence.

For $k\in\{0,\ldots,n-1\}$ with $t_k<T^\e$, $t_k<\tau<\min\{t_{k+1},T^\e\}$, and $\g\in\{0,s_c,\s\}$, set
\begin{align*}
    M_k^\g(\tau) & :=
    \||\nabla|^\g\xwt{R_1}\|_{Y^{\frac{1}{2}}([t_k,\tau])}
    +\||\nabla|^\g R_2\|_{Y^0([t_k,\tau])},                                        \\
    E_k^\g       & :=\|\lna^{\frac{1}{2}}|\nabla|^\g\xwt{R_1}(t_k)\|_{L_x^2(\R^d)}
    +\||\nabla|^\g R_2(t_k)\|_{L_x^2(\R^d)}.
\end{align*}
Let $M_k^\g(t_k) = E_k^\g$. This agrees with the limit of $M_k^\g(\tau)$ as $\tau \to t_k$.

By \eqref{init-R1R2} and the definition of $\xwt{R_1}$, we have
\begin{align} \label{E0}
    E_0^\g=0,\qquad \g\in\{0,s_c,\s\}.
\end{align}

We choose $C\ge1$ sufficiently large to dominate the constants in Proposition~\ref{prop:time-partition}, all coefficients arising from combining \eqref{est-R1} with \eqref{est-R2}, and those arising from combining \eqref{R1-general} with \eqref{R2-general}. We keep $C$ fixed throughout the following argument.

Since $R_1=\Im \xwt{R_1}$, adding \eqref{est-R1} and \eqref{est-R2} on $[t_k,\tau]$ gives
\begin{align} \label{Mk}
    M_k^{s_c}(\tau) \le C E_k^{s_c} + C M_k^{s_c}(\tau)^2 + C \de M_k^{s_c}(\tau) + C C_0^2 \e^2,
\end{align}
for any $k\in\{0,\ldots,n-1\}$.

\medskip

In the sequel of this section, we will give the uniform estimates for $M_{k}^{\g}$ and
show the maximal existence time in Proposition \ref{prop:local-original} $T^\e= \infty$.

\paragraph{First interval: $k=0$.}
Substituting \eqref{E0} into \eqref{Mk}, we obtain
\begin{align} \label{boot-step0}
    M_0^{s_c}(\tau) \le C M_0^{s_c}(\tau)^2 + C \de M_0^{s_c} (\tau) + C C_0^2 \e^2.
\end{align}
Define the bootstrap time by
\begin{align*}
    T_0^*:= \sup \bigl\{ \tau \in  [0,\min\{T^\e,t_{1}\}): M_0^{s_c}(\tau) \le 2 C C_0^2 \e^2 \bigr\}.
\end{align*}
The strict inequality $M_0^{s_c}(0)=0<2 C C_0^2 \e^2$ and continuity give $T_0^*>0$. Since $M_0^{s_c}$ is nondecreasing, the bootstrap bound holds for every $0\le \tau<T_0^*$. Substituting this bound into \eqref{boot-step0}, we obtain
\begin{align} \label{M0-expr}
    M_0^{s_c}(\tau) \le 4 C^3 C_0^4 \e^4 + (2\de C + 1) C C_0^2 \e^2 =   ( 4C^2 C_0^2\epsilon^2 + 2C\delta+1 )C C_0^2\e^2.
\end{align}
Choose $\de$ and $\e_0$ such that
\begin{align} \label{choice:de-0}
    2 \de C = \frac{1}{4},  \qquad 4C^2 C_0^2 \e_0^2 \le \frac{1}{4}.
\end{align}
Then for all $0<\e\le\e_0$, \eqref{M0-expr} gives $M_0^{s_c}(\tau)\le\frac32 C C_0^2\e^2$, which is strictly smaller than $2C C_0^2\e^2$.
Consequently, $T_0^*=\min\{T^\e,t_1\}$ by continuity, and hence
\begin{align} \label{M0-expr-2}
    M_0^{s_c}(\tau)\le 2C C_0^2 \e^2,
    \qquad 0\le\tau<\min\{T^\e,t_1\}.
\end{align}

We next estimate the orders $\g=0,\s$. Using the fixed constant selected before \eqref{Mk}, combining \eqref{R1-general} and \eqref{R2-general}, and using $E_0^\g=0$, we obtain
\begin{align*}
    M_0^\g(\tau) \le C(\de+M_0^{s_c}(\tau))M_0^\g(\tau) + C C_0\e^{\g-s_c}M_0^{s_c}(\tau) + C C_0^2 \e^{\g+2-s_c}, \qquad \g=0,\s.
\end{align*}
Here all factors have already been absorbed into the fixed constant $C$.

Using \eqref{M0-expr-2}, we infer
\begin{align*}
    M_0^\g(\tau) & \le C(\de+2C C_0^2\e^2)M_0^\g(\tau) + 2C^2 C_0^3\e^{\g+2-s_c} + C C_0^2\e^{\g+2-s_c} \\
                 & = C(\de+2C C_0^2\e^2)M_0^\g(\tau) + C C_0^2(1+2C C_0)\e^{\g+2-s_c}, \qquad \g=0,\s.
\end{align*}
By \eqref{choice:de-0}, we have $C\de = \frac18$ and $2C^2 C_0^2\e_0^2\le\frac18$. Hence
\begin{align*}
    C(\de+2C C_0^2\e^2) \le C\de+2C^2 C_0^2\e_0^2 \le \frac14.
\end{align*}
Absorbing the first term on the right-hand side and setting $B_0= 2 C C_0^2(1+2C C_0)$, we obtain
\begin{align}\label{Mg-0}
    M_0^\g(\tau) \le  B_{0} \e^{\g+2-s_c}, \qquad \g=0,\s,\quad 0\le\tau<\min\{T^\e,t_1\}.
\end{align}

If $T^\e<\infty$ and $T^\e\le t_1$, letting $\tau \to T^\e$ in the estimate \eqref{Mg-0} gives
\begin{align*}
    \sum_{\g\in\{0,\s\}}
    \big(
    \||\nabla|^\g\xwt{R_1}\|_{X([0,T^\e))}
    +\||\nabla|^\g R_2\|_{X([0,T^\e))}
    \big)<\infty.
\end{align*}
This contradicts \eqref{continuation-criterion}. Thus, $T^\e>t_1$.

\paragraph{Second interval: $k=1$.}
If $t_1<\infty$, then $T^\e>t_1$. By time continuity and \eqref{M0-expr-2}, we obtain
\begin{align*} %\label{E1}
    E_1^{s_c}\le M_0^{s_c}(t_1)\le 2C C_0^2\e^2.
\end{align*}
Similarly, \eqref{Mg-0} gives
\begin{align}\label{E1-g}
    E_1^\g\le M_0^\g(t_1)\le 2C C_0^2(1+2C C_0)\e^{\g+2-s_c}, \qquad \g=0,\s.
\end{align}

Substituting $E_1^{s_c}\le2C C_0^2\e^2$ into \eqref{Mk}, we obtain
\begin{align}\label{M1}
    M_1^{s_c}(\tau) \le C M_1^{s_c}(\tau)^2 + C\de M_1^{s_c}(\tau) + C(2C+1)C_0^2\e^2, \qquad t_1\le\tau<\min\{T^\e,t_2\}.
\end{align}
Define the bootstrap time by
\begin{align*}
    T_1^* := \sup\bigl\{\tau\in[t_1,\min\{T^\e,t_2\}): M_1^{s_c}(\tau)\le2C(2C+1)C_0^2\e^2\bigr\}.
\end{align*}
Since $M_1^{s_c}(t_1)=E_1^{s_c}\le2C C_0^2\e^2<2C(2C+1)C_0^2\e^2$, continuity gives $T_1^*>t_1$. For $t_1\le\tau<T_1^*$, substituting the bootstrap bound into \eqref{M1} gives
\begin{align*}
    M_1^{s_c}(\tau) & \le 4C^3(2C+1)^2 C_0^4\e^4 + 2C^2\de(2C+1)C_0^2\e^2 + C(2C+1)C_0^2\e^2 \\
                    & = \bigl( 4C^2(2C+1) C_0^2\e^2+2C\de+1 \bigr)C(2C+1)C_0^2\e^2.
\end{align*}
Keeping $\de$ unchanged as chosen in \eqref{choice:de-0}, we decrease $\e_0$ if necessary so that
\begin{align}\label{choice:e-2}
    4C^2(2C+1)C_0^2\e_0^2\le\frac14.
\end{align}
Together with $2C\de = \frac14$, this gives $M_1^{s_c}(\tau)\le\frac32 C(2C+1)C_0^2\e^2$ for $t_1\le\tau<T_1^*$. This is strictly smaller than $2C(2C+1)C_0^2\e^2$.   Consequently, $T_1^*=\min\{T^\e,t_2\}$, and hence
\begin{align}\label{Msc-1}
    M_1^{s_c}(\tau)\le2C(2C+1)C_0^2\e^2, \qquad t_1\le\tau<\min\{T^\e,t_2\}.
\end{align}

We next estimate the orders $\g=0,\s$ at $k=1$. Combining \eqref{R1-general} and \eqref{R2-general}, and using \eqref{Msc-1}, we obtain
\begin{align*}
    M_1^\g(\tau) & \le C E_1^\g + C\bigl(\de+M_1^{s_c}(\tau)\bigr)M_1^\g(\tau) + C C_0\e^{\g-s_c}M_1^{s_c}(\tau) + C C_0^2\e^{\g+2-s_c} \\
                 & \le C E_1^\g + C\bigl(\de+2C(2C+1)C_0^2\e^2\bigr) M_1^\g(\tau) + C C_0^2\bigl(1+2C(2C+1)C_0\bigr)\e^{\g+2-s_c}.
\end{align*}
By \eqref{choice:de-0} and \eqref{choice:e-2}, we have
\begin{align*}
    C\bigl(\de+2C(2C+1)C_0^2\e^2\bigr) \le C\de+2C^2(2C+1)C_0^2\e_0^2 \le \frac18+\frac18=\frac14.
\end{align*}
Absorbing this term and using \eqref{E1-g}, we obtain
\begin{align}
    M_1^\g(\tau) & \le 2C E_1^\g + 2C C_0^2\bigl(1+2C(2C+1)C_0\bigr)\e^{\g+2-s_c} \notag                                                 \\
                 & \le 2C C_0^2\bigl(2C+1+2C(4C+1)C_0\bigr)\e^{\g+2-s_c}, \qquad \g=0,\s,\quad t_1\le\tau<\min\{T^\e,t_2\}. \label{Mg-1}
\end{align}
If $T^\e<\infty$ and $T^\e\le t_2$, letting $\tau\to T^\e$ in \eqref{Mg-1} gives
\begin{align*}
    \sum_{\g\in\{0,\s\}} \bigl(\||\nabla|^\g\xwt{R_1}\|_{X([t_1,T^\e))} + \||\nabla|^\g R_2\|_{X([t_1,T^\e))}\bigr) < \infty.
\end{align*}
Together with the bounds on $[0,t_1]$, this contradicts \eqref{continuation-criterion}. Thus, $T^\e>t_2$.

\paragraph{Remaining intervals: $2\le k\le n-1$.}
%If $t_2<\infty$, we continue the argument by induction. 

By induction, we assume for $ j = 0, 1, \cdots, k-1$ that
\ba\label{Mk-induction-ass}
M_j^{s_{c}}(\tau)  \leq A_{j} \e^{2}, \quad    M_j^{\gamma}(\tau)  \leq B_{j} \e^{\g+2-s_{c}}, \quad \g = 0, \sigma.
\ea
We have shown that \eqref{Mk-induction-ass} holds for $j=0,1$ with
\begin{align*}
    A_0=2C C_0^2, \qquad         & B_0=2C C_0^2 (1 + 2C C_0),                   \\
    A_1=2C C_0^2 (2C+1),  \qquad & B_1=2C C_0^2 \bigl(2C+1 + 2C(4C+1)C_0\bigr).
\end{align*}
By time continuity and the estimates on $I_{k-1}$, we have
\begin{align*}
    E_k^{s_c} \le A_{k-1}\e^2, \qquad E_k^\g \le B_{k-1}\e^{\g+2-s_c}, \quad \g=0,\s.
\end{align*}
Substituting the first bound into \eqref{Mk}, we obtain
\begin{align*}
    M_k^{s_c}(\tau) \le C M_k^{s_c}(\tau)^2 + C\de M_k^{s_c}(\tau) + (C A_{k-1}+C C_0^2)\e^2.
\end{align*}
Set
\begin{align*}
    T_k^* := \sup \bigl\{ \tau \in [t_k,\min\{T^\e,t_{k+1}\}) : M_k^{s_c}(\tau) \le A_k\e^2 \bigr\}.
\end{align*}
Here we choose
\begin{align*} %\label{cond:Ak}
    A_k=2C A_{k-1}+ 2C C_0^2.
\end{align*}
The bootstrap bound and the definition of $A_k$ yield
\begin{align*}
    M_k^{s_c}(\tau) & \le C A_k^2\e^4 + C\de A_k\e^2 + \frac12 A_k\e^2 = \bigl(C A_k\e^2 + C\de + \tfrac12\bigr) A_k\e^2.
\end{align*}
Under the condition
\begin{align} \label{cond:ek}
    C A_k\e_0^2 \le \frac18,
\end{align}
and recalling that $C\de \le \frac18$, we obtain $M_k^{s_c}(\tau) \le \frac34 A_k\e^2$. The same continuity argument as for $k=0,1$ therefore gives $T_k^*=\min\{T^\e,t_{k+1}\}$, and hence
\begin{align*}
    M_k^{s_c}(\tau) \le A_k\e^2, \qquad t_k \le \tau < \min\{T^\e,t_{k+1}\}.
\end{align*}

At orders $\g=0,\s$, combining \eqref{R1-general} and \eqref{R2-general} with the profile bounds and the preceding critical-order estimate gives
\begin{align*}
    M_k^\g(\tau) \le C E_k^\g + C(\de+A_k\e^2)M_k^\g(\tau) + C C_0 A_k\e^{\g+2-s_c} + C C_0^2\e^{\g+2-s_c}.
\end{align*}
Since $C(\de+A_k\e^2)\le\frac14$, absorbing this term and using the induction hypothesis, we obtain
\begin{align*}
    M_k^\g(\tau) & \le 2C E_k^\g + (2C C_0 A_k+2C C_0^2)\e^{\g+2-s_c}  \le  B_{k} \e^{\g+2-s_c}, \quad \g=0,\s,
\end{align*}
where we set
\begin{align}\label{cond:Bk}
    B_k = 2C B_{k-1}+2C C_0 A_k+2C C_0^2.
\end{align}
The bounds of $M_{j}^{\gamma}$ for $j=0,\cdots, k$ contradict \eqref{continuation-criterion} if $T^\e<\infty$ and $T^\e\le t_{k+1}$. Thus $T^\e>t_{k+1}$ when $t_{k+1}<\infty$. At this finite endpoint, time continuity gives
\begin{align*}
    E_{k+1}^{s_c} \le A_k\e^2, \qquad E_{k+1}^\g \le B_k\e^{\g+2-s_c}, \quad \g=0,\s,
\end{align*}
which supplies the initial bounds for the next interval. If $t_{k+1}=\infty$, the same contradiction gives $T^\e=\infty$.

We now estimate $A_k$. Since $A_0=2C C_0^2$ and $A_k=2C A_{k-1}+2C C_0^2$, by induction we have
\begin{align*}
    A_k=2C C_0^2\sum_{j=0}^k(2C)^j \le 2C C_0^2(2C+1)^k.
\end{align*}
With the choice of $\de$ in \eqref{choice:de-0}, Proposition~\ref{prop:time-partition} gives
\begin{align} \label{bound:n}
    n & \le 2 + \lfloor 2 \sum_{j=1}^2 (64 C^2 C_0)^{p_j} \rfloor  \le 2 + 2 \sum_{j=1}^2 (64 C^2 C_0)^{p_j}     \le 4 (64 C^2)^{\frac{d+2}{2}} \bigl(1 + C_0^{\frac{d+2}{2}}\bigr).
\end{align}
%Here we used $C\ge1$, $p_2\le p_1=\frac{d+2}{2}$, and $C_0^{p_2}\le 1+C_0^{p_1}$.

To ensure \eqref{cond:ek} on every interval $I_k$, we may further decrease $\e_0$ if necessary so that
\begin{align*} %\label{final-e0}
    2C^2 C_0^2(2C+1)^{\lfloor N_\de \rfloor}\e_0^2\le\frac18.
\end{align*}
Then
\begin{align*}
    C A_k\e_0^2 \le 2C^2 C_0^2(2C+1)^k\e_0^2 \le\frac18, \qquad 0\le k\le n-1 \le N_\de.
\end{align*}
The bound for $n$ in Proposition~\ref{prop:time-partition} also makes this choice of $\e_0$ uniform when $C_0$ is bounded above. The preceding induction therefore applies to all intervals and gives $T^\e=\infty$. Consequently,
\begin{align*}
    M_k^{s_c}(\tau)\le2C C_0^2(2C+1)^k\e^2, \qquad \forall \, t_k\le\tau<t_{k+1},\quad \forall\,  0\le k\le n-1.
\end{align*}

We next estimate $B_k$. Iterating \eqref{cond:Bk} and using $B_0=2C C_0^2+2C C_0 A_0$, we obtain
\begin{align*}
    B_k=2C C_0^2\sum_{j=0}^k(2C)^{k-j}+2C C_0\sum_{j=0}^k(2C)^{k-j}A_j.
\end{align*}
Substituting the expression for $A_j$ gives
\ba
B_k  =2C C_0^2\sum_{j=0}^k(2C)^j+4C^2 C_0^3\sum_{j=0}^k\sum_{l=0}^{k-j}(2C)^{j+l} \leq  2C C_0^2\bigl(1+2C C_0(2C+1)\bigr)(2C+1)^k.
\nn \ea
Thus,
\begin{align*}
    M_k^\g(\tau)\le2C C_0^2\bigl(1+2C C_0(2C+1)\bigr)(2C+1)^k\e^{\g+2-s_c}, \qquad \g=0,\s,
\end{align*}
for all $t_k\le\tau<t_{k+1}$ and all $0\le k\le n-1$.

\subsection{End of the proof}

Letting $\tau\to t_{k+1}$ in the interval estimates, with $\tau\to \infty$ on the last interval, gives the corresponding bounds on $I_k$. Summing over the intervals, we obtain
\begin{align*}
    \||\nabla|^\g\xwt{R_1}\|_{Y^{\frac12}(\R_+)} + \||\nabla|^\g R_2\|_{Y^0(\R_+)} & \le \sum_{k=0}^{n-1}\bigl(\||\nabla|^\g\xwt{R_1}\|_{Y^{\frac12}(I_k)} + \||\nabla|^\g R_2\|_{Y^0(I_k)}\bigr) \\
                                                                                   & \le 2C C_0^2\bigl(1+2C C_0(2C+1)\bigr)\e^{\g+2-s_c}\sum_{k=0}^{n-1}(2C+1)^k                                  \\
                                                                                   & \le 2C C_0^2\bigl(1+2C C_0(2C+1)\bigr)(2C+2)^n\e^{\g+2-s_c}, \qquad \g=0,\s.
\end{align*}

By \eqref{local-class}, the global bounds above, and the identities in \eqref{relRr2}, we have
\begin{align*}
    R_1\in C_b(\R_+;H^{\s+\frac12}(\R^d)), \quad \pa_t R_1\in C_b(\R_+;H^{\s-\frac12}(\R^d)), \quad R_2\in C_b(\R_+;H^\s(\R^d)).
\end{align*}
Returning to the original variables by \eqref{scale} and adding the profiles, we obtain
\begin{align*}
    \phi\in C_b(\R_+;H^{\s+\frac12}(\R^d)), \quad \pa_t\phi\in C_b(\R_+;H^{\s-\frac12}(\R^d)), \quad \psi\in C_b(\R_+;H^\s(\R^d)).
\end{align*}
For the approximation error, using \eqref{scale}, the identity $R_1=\Im\xwt{R_1}$, and the global remainder bounds, we obtain
\begin{align*}
    \||\nabla|^\g r_1\|_{L_t^\infty L_x^2(\R_+)} + \||\nabla|^\g r_2\|_{L_t^\infty L_x^2(\R_+)} & \le \e^{s_c-\g}\bigl(\||\nabla|^\g\xwt{R_1}\|_{Y^{\frac12}(\R_+)} + \||\nabla|^\g R_2\|_{Y^0(\R_+)}\bigr) \\
                                                                                                & \le 2C C_0^2\bigl(1+2C C_0(2C+1)\bigr)(2C+2)^n\e^2, \qquad \g=0,\s.
\end{align*}
Therefore,
\begin{align*} %\label{est:rl}
    \|\phi(t)-e^{\frac{it}{\e^2}}\phi_s(t)
    -e^{-\frac{it}{\e^2}}\xwo{\phi_s}(t)\|_{H^\s(\R^d)}
    +\|\psi(t)-\psi_s(t)\|_{H^\s(\R^d)}
    \le \xwt{C_0} C_0^2 \e^2,\qquad t\ge 0,
\end{align*}
where, by \eqref{bound:n}, we may take
\begin{align*}
    \xwt{C_0} := 2 C  \bigl(1+2C C_0(2C+1)\bigr)(2C+2)^{4(64 C^2)^{\frac{d+2}{2}}\bigl(1+C_0^{\frac{d+2}{2}}\bigr)}.
\end{align*}
This constant depends only on $d,\s,\l,C_0$. This completes the proof of Theorem~\ref{thm:1}.

\medskip

We now prove Corollary~\ref{cor:general-Phi1} by comparing the Schr\"odinger-wave profiles with initial velocities $0$ and $\Phi_1$.
\begin{proof}[Proof of Corollary~\ref{cor:general-Phi1}]
    The profiles $\phi_s^{\Phi_1}$ and $\phi_s$ satisfy the same linear Schr\"odinger-wave equation, and their difference has initial data
    \begin{align*}
        (\phi_s^{\Phi_1}-\phi_s)(0)=0,
        \qquad
        \pa_t(\phi_s^{\Phi_1}-\phi_s)(0)=\Phi_1.
    \end{align*}
    Multiplication by $e^{\frac{it}{\e^2}}$ transforms the equation for their difference into the homogeneous Klein--Gordon equation:
    \begin{align*}
        \left(\e^2\pa_t^2-\Delta+\e^{-2}\right)
        [e^{\frac{it}{\e^2}}
            (\phi_s^{\Phi_1}-\phi_s) ]=0.
    \end{align*}
    Since the difference vanishes at $t=0$, the transformed function also has zero initial displacement and initial velocity $\Phi_1$. The Klein--Gordon solution formula therefore gives
    \begin{align*}
         & e^{\frac{it}{\e^2}}(\phi_s^{\Phi_1}-\phi_s)(t) = \e^2(1-\e^2\Delta)^{-\frac12}
        \sin \big(\frac{t}{\e^2}(1-\e^2\Delta)^{\frac12}\big)\Phi_1.
    \end{align*}
    By Plancherel's theorem,
    \begin{align*}
         & \sup_{t\ge 0}
        \|2\Re\big(e^{\frac{it}{\e^2}}
        (\phi_s^{\Phi_1}-\phi_s)(t)\big) \|_{H^\s(\R^d)} \le 2\e^2\|\Re\Phi_1\|_{H^\s(\R^d)}
        \le 2 \e^2 \|\Phi_1\|_{H^\s(\R^d)}.
    \end{align*}
    Combining this estimate with \eqref{eq:main-1} and the triangle inequality, we obtain, for every $t\ge 0$ and $0<\e\le\e_0$,
    \begin{align*}
         & \|\phi(t)-e^{\frac{it}{\e^2}}\phi_s^{\Phi_1}(t)
        -e^{-\frac{it}{\e^2}}\xwo{\phi_s^{\Phi_1}}(t) \|_{H^\s(\R^d)}
        +\|\psi(t)-\psi_s(t)\|_{H^\s(\R^d)}                                       \\
         & \quad\le  \xwt{C_0} C_0^2 \e^2+
        \|2\Re \big(e^{\frac{it}{\e^2}}
        (\phi_s^{\Phi_1}-\phi_s)(t)\big) \|_{H^\s(\R^d)}                          \\
         & \quad \le \big( \xwt{C_0} C_0^2 +2\|\Re\Phi_1\|_{H^\s(\R^d)}\big)\e^2.
    \end{align*}
    This proves Corollary \ref{cor:general-Phi1} with $C_{\Phi_1}=\xwt{C_0}C_0^2+2\|\Re\Phi_1\|_{H^\s(\R^d)}$.
\end{proof}

\section*{Acknowledgements}
Y. L. and Z. Z. are partially supported by Basic Research Program of Jiangsu (Grant No. BK20240058) and by National Natural Science Foundation of China (Grant No. 12671257). Y. L. is also partially supported by the Fundamental Research Funds for the Central Universities (Grant No. 2026300394).

\section*{Conflict of interest}
The authors declare no conflict of interest in this paper.

\end{document}